# RECURSIVE MONTE CARLO FILTERS: ALGORITHMS AND THEORETICAL ANALYSIS

By Hans R. Künsch

*ETH Zürich*

Recursive Monte Carlo filters, also called particle filters, are a powerful tool to perform computations in general state space models. We discuss and compare the accept–reject version with the more common sampling importance resampling version of the algorithm. In particular, we show how auxiliary variable methods and stratification can be used in the accept–reject version, and we compare different resampling techniques. In a second part, we show laws of large numbers and a central limit theorem for these Monte Carlo filters by simple induction arguments that need only weak conditions. We also show that, under stronger conditions, the required sample size is independent of the length of the observed series.

## 1. State space and hidden Markov models.

A *general state space* or *hidden Markov model* consists of an unobserved state sequence $(X_t)$ and an observation sequence $(Y_t)$ with the following properties:

*State evolution*: $X_0, X_1, X_2, \ldots$ is a Markov chain with $X_0 \sim a_0(x)\, d\mu(x)$ and

$$X_t | X_{t-1} = x_{t-1} \sim a_t(x_{t-1}, x)\, d\mu(x).$$

*Generation of observations*: Conditionally on $(X_t)$, the $Y_t$'s are independent and $Y_t$ depends on $X_t$ only with

$$Y_t | X_t = x_t \sim b_t(x_t, y)\, d\nu(y).$$

These models occur in a variety of applications. Linear state space models are equivalent to ARMA models (see, e.g., [16]) and have become popular under the name of structural models (see, e.g., [17]). Nonlinear state space models occur in finance (stochastic volatility; see, e.g., [27]), in various fields of engineering (speech, tracking and control problems; see, e.g., [12]), in biology (ion channels, DNA and protein sequences) and in geophysics (rainfall at a network of stations, data assimilation). A more detailed survey with many references is given in [20].

In order to apply these models, two kinds of problems have to be solved: Inference about the states based on a stretch of observed values $y_{s:t} = (y_u, s \leq u \leq t)$ for a given model, that is, $a_t$ and $b_t$ known (this is called prediction, filtering and smoothing), and inference about unknown parameters in $a_t$, $b_t$. From a statistical point of view, the latter problem is maybe of greater interest, but fast and reliable algorithms for the former are a prerequisite for computing maximum likelihood or Bayesian estimators. The reason for this is briefly mentioned in Section 2.1. This paper is therefore entirely devoted to algorithms for filtering, prediction and smoothing.

Section 2 recalls the basic recursions for filtering, prediction and smoothing. Section 3 discusses algorithmic aspects of sequential Monte Carlo methods to implement these recursions. Most algorithms in the literature, beginning with the pioneering paper by Gordon, Salmond and Smith [15], use the sampling importance resampling idea of Rubin [26]. An exception is Hürzeler and Künsch [18] who use the accept–reject method instead. Here we show how some ideas like stratification and an auxiliary variable method of Pitt and Shephard [23] can be adapted to rejection sampling, and we give new results on the performance of systematic resampling methods. In addition, we hope that our view of classifying and comparing approaches is useful.

Section 4 presents results on the convergence of the method as the number of Monte Carlo replicates tends to infinity. We discuss both laws of large numbers and a central limit theorem. Recently, many similar results have been published; see, for example, [4, 8, 21]. The distinctive features of our presentation here are the weakness of conditions, the use of the total variation distance to measure the difference between the approximate and the true filter density and the simplicity of the techniques used. We basically show that most results follow by induction, in accordance with the recursive nature of the algorithm. The complications that occur are due to a counterintuitive property of Bayes' formula; see Lemma 3.6(ii) in [20]. As a consequence, although one can obtain consistency with very few conditions on the model, the required sample size seems to grow exponentially with the number of time steps. For results that guarantee that the required sample size is independent of the number of time steps (or grows at most logarithmically), one has to use induction over several time steps which requires rather strong conditions on the dynamics of the states. At the end, we give some results for the case where $(X_t)$ is a continuous time process and the sampling rate of the observations increases.

**2. Filtering and smoothing recursions.** In general, we will use the symbol $p$ as the generic notation for a conditional density of its arguments. However, for the conditional density of $X_t$ given $Y_{1:s} = y_{1:s}$, we use the notation $f_{t|s}(x_t|y_{1:s})$. The three cases $s < t$, $s = t$ and $s > t$ are called prediction, filtering and smoothing, respectively.

The dependence structure of a state space model can be represented by the following directed acyclic graph:

$$
\begin{array}{ccccccccc}
\cdots & \longrightarrow & X_{t-1} & \longrightarrow & X_t & \longrightarrow & X_{t+1} & \longrightarrow & \cdots \\
 & & \downarrow & & \downarrow & & \downarrow & & \\
\cdots & & Y_{t-1} & & Y_t & & Y_{t+1} & & \cdots
\end{array}
$$

From this, various conditional independence properties follow which are used together with the law of total probability and Bayes theorem to derive recursions for the filter, prediction and smoothing densities. These are well known; see, for example, [20], Section 3.3, and we state them without proofs.

The most important result is the following recursion for the filter density:

*Propagation.* From the filter density, we obtain the one-step-ahead prediction density:

$$
f_{t|t-1}(x_t|y_{1:t-1}) = \int f_{t-1|t-1}(x|y_{1:t-1}) a_t(x, x_t)\, d\mu(x). \tag{1}
$$

*Update.* From the-one-step ahead prediction density, we obtain the filter density one time step later:

$$
\begin{aligned}
f_{t|t}(x_t|y_{1:t}) &= \frac{f_{t|t-1}(x_t|y_{1:t-1}) b_t(x_t, y_t)}{\int f_{t|t-1}(x|y_{1:t-1}) b_t(x, y_t)\, d\mu(x)} \\
&\propto f_{t|t-1}(x_t|y_{1:t-1}) b_t(x_t, y_t).
\end{aligned} \tag{2}
$$

In parts of the literature, for example, in [8], $Y_t$ depends on $X_{t-1}$ and not on $X_t$. Then the filter density is, in our setup, the prediction density which should be kept in mind when comparing formulae.

2.1. *Prediction of observations and likelihood.* The denominator in the update step (2) is the conditional density of $Y_t$ given $Y_{1:t-1}$:

$$
p(y_t|y_{1:t-1}) = \int f_{t|t-1}(x|y_{1:t-1}) b_t(x, y_t)\, d\mu(x). \tag{3}
$$

If $f_{t|t-1}$ is available, we thus can obtain the likelihood from

$$
p(y_{1:T}) = \prod_{t=1}^{T} p(y_t|y_{1:t-1}).
$$

A different representation of the likelihood is obtained by marginalization:

$$p(y_{1:T}) = \int a_0(x_0) \prod_{t=1}^{T} a_t(x_{t-1}, x_t) b_t(x_t, y_t) \prod_{t=0}^{T} d\mu(x_t).$$

From this, the likelihood ratio can be expressed as an expectation with respect to the smoothing distribution; see, for example, [19].

2.2. *Smoothing.* The filter densities can also be used for the smoothing problem since conditional on $y_{1:T}$, $(X_T, X_{T-1}, \ldots, X_0)$ is an inhomogeneous Markov chain with starting density $f_{T|T}$ and backward transition densities

$$p(x_t | x_{t+1}, y_{1:T}) = p(x_t | x_{t+1}, y_{1:t}) \propto a_{t+1}(x_t, x_{t+1}) f_{t|t}(x_t | y_{1:t}). \tag{4}$$

This is also the basis for the forward-filtering–backward-sampling algorithm; see [13], equation (20). From (4), we can derive, in particular, a backward recursion for $f_{t|T}$.

2.3. *Recursive filtering in operator notation.* A compact notation for the filter recursion which will be useful later on is

$$f_{t|t}(\cdot | y_{1:t}) = B(A_t^* f_{t-1|t-1}(\cdot | y_{1:t-1}), b_t(\cdot, y_t)). \tag{5}$$

Here

$$A_t^* f(x) = \int f(x') a_t(x', x) \, d\mu(x')$$

is the Markov transition operator, and

$$B(f, b)(x) = \frac{f(x) b(x)}{\int f(x) b(x) \, d\mu(x)}$$

is the Bayes operator that assigns the posterior to a prior $f$ and a likelihood $b$. The operators $A_t^*$ and $B(\cdot, b)$ map the space of densities into itself, but they can be extended to the space of probability distributions.

2.4. *Implementation of recursions.* If $X_t$ is discrete with $M$ possible values, integrals are sums and the recursions need $O(TM^2)$ operations. In a linear Gaussian state space model, all $f_{t|s}$ are Gaussian, and their means and variances are computed with the Kalman filter and smoother.

In practically all other cases, the recursions are difficult to compute. Analytical approximations like the extended Kalman filter are not satisfactory, and numerical integration is problematic in high dimensions. Much current interest focuses on Monte Carlo methods. Standard Markov chain Monte Carlo can be used, but it lacks a recursive implementation. There has been considerable interest in recursive Monte Carlo methods in recent years; see, for example, [12].

**3. Algorithms for recursive Monte Carlo filtering.** The following is the key observation: $A_t^* f$ is difficult to compute, but easy to sample from if we can sample from $f$ and $a_t(x,\cdot)$. This allows us to generate recursively a sequence of samples ("particles") $(x_{j,t}; j=1,\dots,N, t=0,1,\dots)$ with approximate distribution $f_{t|t}$ as follows: If $(x_{j,t-1})$ is available, we can replace

$$A_t^* f_{t-1|t-1}(x_t|y_{1:t-1}) = \int f_{t-1|t-1}(x|y_{1:t-1}) a_t(x,x_t)\,d\mu(x)$$

by

$$\frac{1}{N}\sum_{j=1}^{N} a_t(x_{j,t-1},x_t).$$

Therefore, we sample $(x_{j,t})$ from the distribution with density

$$f_{t|t}^N(\cdot|y_{1:t}) \propto b_t(\cdot,y_t)\frac{1}{N}\sum_{j=1}^{N} a_t(x_{j,t-1},\cdot). \tag{6}$$

In this section we discuss methods to sample from such a density. We simplify the notation somewhat and write the target density as

$$f^N(x) \propto f_u^N(x) = b(x)\sum_{j=1}^{N} a(j,x) \tag{7}$$

(subscript $u$ for unnormalized). We will call $b$ the likelihood and $N^{-1}\sum_j a(j,x)$ the prior. In the filtering context, the prior is the approximate prediction density. For later use, we also introduce

$$\beta_j = \int a(j,x)b(x)\,d\mu(x),$$

which is in the filtering context equal to the conditional density of $Y_t$ given $X_{t-1} = x_{j,t-1}$. We assume that we have good methods to generate samples from $a(j,\cdot)$ for any $j$. The methods we discuss fall into two categories: accept–reject and importance sampling with an additional resampling step.

3.1. *Accept–reject methods.* The accept–reject method for sampling from the density (7) produces values $X$ according to a proposal $\rho$, and if $X=x$ accepts it with probability

$$\pi(x) = \frac{f_u^N(x)}{M\rho(x)}. \tag{8}$$

Here $M$ is an upper bound for the ratio $f_u^N(x)/\rho(x)$:

$$M \geq \sup_x \frac{f_u^N(x)}{\rho(x)}.$$

The most obvious proposal $\rho(x)$ is the prior, that is,

$$\rho(x) = \frac{1}{N}\sum_{j=1}^{N} a(j,x). \tag{9}$$

Then the evaluation of the acceptance probabilities $\pi(x)$ is easy as long as $b$ is bounded. In order to sample from (9), we first choose an index $J$ uniformly from $\{1,\dots,N\}$, and given $J=j$, we sample $X$ from $a(j,x)$. Note that, in this case, the densities $a(j,x)$ need not be available in analytic form; we only have to be able to sample from them. This is of interest in discretely observed diffusion models.

The average acceptance probability of this algorithm is $\int \rho(x)\pi(x)\,d\mu(x) = \sum_j \beta_j/M$. In particular, if $\rho$ is the prior and if we use the smallest value of $M$, it is equal to

$$\frac{\sum_{j=1}^{N}\beta_j}{N\sup_x b(x)}.$$

This is low if the likelihood is more informative (concentrated) than the prior, or if the likelihood and the prior are in conflict. We discuss here some modifications and tricks that can alleviate this problem in some situations.

3.1.1. *The mixture index as auxiliary variable.* Other proposal distributions than the prediction density can, of course, lead to higher acceptance rates, but usually it is difficult to compute a good upper bound $M$, and the evaluation of the acceptance probability $\pi(x)$ is complicated due to the sum over $j$. A way to avoid at least the last problem is based on an idea by Pitt and Shephard [23]. Namely, we can generate first an index $J$ according to a distribution $(\tau_j)$ and given $J=j$, a variable $X$ according to a density $\rho(j,x)$. We then accept the generated pair $(j,x)$ with probability

$$\pi(j,x) = \frac{a(j,x)b(x)}{M\tau_j\rho(j,x)}, \tag{10}$$

where now

$$M \geq \sup_{j,x}\frac{a(j,x)b(x)}{\tau_j\rho(j,x)}.$$

If the pair is accepted, we simply discard $j$ and keep $x$, and otherwise we generate a new pair. Because the accepted pairs $(J,X)$ have distribution

$$\frac{a(j,x)b(x)}{\sum_j\beta_j},$$

the marginal distribution of $X$ is the target (7). If we take $\tau_j = 1/N$ and $\rho(j,x) = a(j,x)$, we obtain the usual algorithm discussed before, but one will try to increase the acceptance rate by other choices.

Because $j$ runs over a finite set, we will usually take

$$M = \max_j \frac{M_j}{\tau_j}, \qquad \text{where } M_j \geq \sup_x \frac{a(j,x)b(x)}{\rho(j,x)}.$$

LEMMA 1. *For a given choice of densities $\rho(j,x)$ and bounds $M_j$, the average acceptance probability is less than or equal to $\sum \beta_j / \sum M_j$, with equality iff $\tau_j \propto M_j$.*

PROOF. The average acceptance probability is

$$\sum_j \int \pi(j,x)\tau_j \rho(j,x)\mu(dx) = \frac{1}{M}\sum_j \beta_j = \sum_j \beta_j \left(\max_k \frac{M_k}{\tau_k}\right)^{-1}.$$

Clearly,

$$\max_k \frac{M_k}{\tau_k} = \sum_j \tau_j \max_k \frac{M_k}{\tau_k} \geq \sum_j M_j,$$

with equality iff $M_k/\tau_k$ is constant. □

If $\rho(j,x) = a(j,x)$, the optimal $\tau_j$'s are thus constant. This is somewhat surprising since one could conjecture that it is better to give higher probability to those indices $j$ for which the mass of $a(j,x)$ is close to $\arg\sup b(x)$.

The crucial point in implementing this algorithm is the choice of the densities $\rho(j,\cdot)$. Lemma 1 implies that, for a high acceptance probability, all $M_j$'s should be small, that is, each $\rho(j,x)$ should be a good proposal distribution for the density $a(j,x)b(x)/\beta_j$. Ideally, we would choose that density itself. But then $M_j$ must be close to the normalizing constant $\beta_j$ which typically is not available in closed form. A more practical approach chooses a parametric family $(\rho(\theta,x))$, where we have available tight upper bounds

$$M(j,\theta) \geq \sup_x \frac{a(j,x)b(x)}{\rho(\theta,x)}.$$

We then optimize over $\theta$, that is,

$$\rho(j,x) = \rho(\theta_j,x), \qquad \text{where } \theta_j \approx \arg\min_\theta M(j,\theta).$$

Note that it is not necessary to find the optimal $\theta$ exactly, but $M(j,\theta)$ should be a true upper bound. By choosing the family $(\rho(\theta,x))$ such that it contains all densities $a(j,x)$, we can make sure that the acceptance probability is at least as high as with the usual algorithm.

EXAMPLE. The simplest stochastic volatility model, see, for example, [27], is obtained if we take for $a(j,\cdot)$ the normal density with mean $m_j$ and variance $\sigma^2$ and for $b$ the likelihood of a $\mathcal{N}(0,\exp(x))$ random variable $Y$,

$$b(x) = b(x,y) = \exp\left(-\frac{x}{2} - \frac{y^2}{2}\exp(-x)\right).$$

If we choose as $\rho(\theta,\cdot)$ the normal density with mean $\theta$ and variance $\sigma^2$, we can compute the supremum of

$$\log\frac{a(j,x)b(x)}{\rho(\theta,x)} = -\frac{x}{2} - \frac{\theta - m_j}{\sigma^2}x - \frac{y^2}{2}\exp(-x) - \frac{m_j^2 - \theta^2}{2\sigma^2}$$

over $x$. It is equal to

$$\frac{\sigma^2}{2}\delta^2 + m_j\delta - \left(\frac{1}{2} + \delta\right)(1 + \log y^2) + \left(\frac{1}{2} + \delta\right)\log(1 + 2\delta),$$

provided $\delta = (\theta - m_j)/\sigma^2 \geq -1/2$ (otherwise the function is unbounded above). Minimizing this expression with respect to $\delta$ subject to $\delta \geq -1/2$ leads to a nonlinear equation which has no closed form solution. Using $\log(1 + 2\delta) \leq 2\delta$, we obtain a quadratic upper bound which is minimized by

$$\theta_j = m_j + \frac{\sigma^2}{2}\max\left(-1, \frac{2}{4 + \sigma^2}(\log y^2 - m_j)\right).$$

This choice of $\theta_j$ may be slightly suboptimal, but because the bound is sharp for $\theta = m_j$, that is, $\delta = 0$, we still can guarantee a higher acceptance probability than with the usual method. In practice, the gain can be dramatic if $|y|$ is small.

The above choice of $\theta_j$ is somewhat different from the suggestion

$$\theta_j = m_j + \frac{\sigma^2}{2}(y^2\exp(-m_j) - 1)$$

in [27], page 285. In addition, the choices for $\tau_j$ differ.

3.1.2. *Balanced sampling.* Besides reducing the acceptance rate, we can also try to reduce the variance by using a more balanced sampling: The target $f^N$ is a mixture of $N$ components, and the variance is reduced if the different components in the mixture are represented with the correct proportions. This idea has received much attention in the sampling importance sampling context; see Section 3.2.1 below and the references given there. We have not seen this idea in the accept–reject context. Consider the estimation of

$$m(\psi) = \int f^N(x)\psi(x)\,d\mu(x) = \frac{\sum_j \beta_j m_j(\psi)}{\sum_j \beta_j},$$

where

$$m_j(\psi) = \frac{\int \psi(x) a(j,x) b(x)\, d\mu(x)}{\beta_j}$$

and $\psi$ is a bounded "test function." If $(X_i)$ is an i.i.d. sample from $f^N$, the estimator

$$\hat{m}(\psi) = \frac{\sum_{j=1}^N \psi(X_j)}{N}$$

has variance

$$\frac{1}{N}\sigma^2(\psi) = \frac{1}{N} \frac{\sum_j \int (\psi(x) - m(\psi))^2 a(j,x) b(x)\, d\mu(x)}{\sum_j \beta_j}.$$

A method to reduce this variance replaces the random selection of an index $J$ by a more systematic procedure. Namely, we can propose simultaneously $N$ values, one each from the density $a(j,x)$, and decide whether to accept each of them independently. We repeat the procedure until the total of accepted values is at least $N$. If we need exactly $N$ values, we can select them at random. We therefore consider the estimator

$$\tilde{m}(\psi) = \frac{\sum_{i=1}^T \sum_{j=1}^N \psi(X_{ij}) \mathbb{1}_{[U_{ij} < b(X_{ij})]}}{\sum_{i=1}^T \sum_{j=1}^N \mathbb{1}_{[U_{ij} < b(X_{ij})]}},$$

where $(X_{ij}, U_{ij}; 1 \le j \le N, i = 1, 2, \dots)$ are independent random variables with $X_{ij} \sim a(j, \cdot)$, $U_{ij}$ uniform on $(0, \sup b(x))$, and $T$ is the smallest integer such that the denominator is at least $N$.

In order to compute the variance of $\tilde{m}(\psi)$ approximately, we use

$$\tilde{m}(\psi) - m(\psi) = \frac{\sum_{i=1}^T \sum_{j=1}^N (\psi(X_{ij}) - m(\psi)) \mathbb{1}_{[U_{ij} < b(X_{ij})]}}{\sum_{i=1}^T \sum_{j=1}^N \mathbb{1}_{[U_{ij} < b(X_{ij})]}}.$$

For simplicity, we assume that $\sup b(x) = 1$. Then, by Wald's identity the denominator has expected value

$$\mathbf{E}[T] \sum_{j=1}^N \beta_j.$$

In particular, the expected number of random variables that have to be generated is essentially the same as with basic i.i.d. rejection sampling. Similarly, the numerator has expectation zero and variance

$$\begin{aligned}
&\mathbf{E}(T) \sum_{j=1}^N \mathrm{Var}((\psi(X_{1j}) - m(\psi)) \mathbb{1}_{[U_{1j} < b(X_{1j})]}) \\
&\quad = \mathbf{E}(T) \left( \sigma^2(\psi) \sum_j \beta_j - \sum_j \beta_j^2 (m_j(\psi) - m(\psi))^2 \right).
\end{aligned}$$

Assuming the denominator to be approximately constant and equal to $N$ (which is reasonable if the expected number of accepted values in each round of proposals is small), we obtain the approximation

$$\mathbf{E}[\hat{m}(\psi)] \approx m(\psi), \qquad \mathrm{Var}(\hat{m}(\psi)) \approx \frac{1}{N}\left(\sigma^2(\psi) - \frac{\sum_j \beta_j^2 (m_j(\psi) - m(\psi))^2}{\sum \beta_j}\right).$$

The second term thus quantifies the gain of the method.

3.2. *Sampling importance resampling.* This method generates $(z_k; 1 \leq k \leq R)$ according to some proposal $\rho$ and selects from these a sample of size $N$ with inclusion probabilities

$$\pi(z_k) \propto \frac{b(z_k) \sum_{j=1}^N a(j, z_k)}{\rho(z_k)}. \tag{11}$$

The resampling need not be made at random. We will discuss below alternative methods with reduced variability. The standard proposal is again the prior (9), leading to the original proposal in [15].

This method has difficulties if the sampling probabilities $\pi(z_k)$ are very unbalanced since this leads to many ties in the final sample. Typically, this occurs in situations where the prior and the likelihood are in conflict, that is, when the acceptance rate in rejection sampling is low. Choosing $R$ much bigger than $N$ reduces the number of ties, but at the expense of longer computations. Note that rejection sampling is an automatic way of choosing $R$ such that all ties are avoided. There are also possibilities for reusing the rejected variables for estimating the current filter distribution more accurately; see Section 3.3.3 of [25].

Most of the ideas discussed in connection with rejection sampling can also be used here. The idea of Pitt and Shephard [23] to include explicitly an index $J$ was originally developed for this case. It proposes a sample $(j_k, z_k)$ of size $R$ with distribution $\tau_j \rho(j, x)$ and then selects a sample of size $N$ with inclusion probabilities

$$\pi(z_k, j_k) \propto \frac{b(z_k) a(j_k, z_k)}{\tau_{j_k} \rho(j_k, z_k)}.$$

In contrast to rejection sampling, combining $\rho(j, \cdot) = a(j, \cdot)$ with unequal $\tau_j$'s is a promising idea here. For instance, we can take $\tau_j$ to be proportional to $b(m_j)$, where $m_j$ is the mean or the median of $a(j, \cdot)$. If all $a(j, \cdot)$'s have a small spread (relative to the scale at which $b$ varies), then most $\pi(z_k, j_k)$'s will be approximately equal and therefore $R = N$ is sufficient.

3.2.1. *Balanced sampling.* Both in the proposal and in the resampling step, we have to select indices from a given distribution. In the former case, this distribution is $(\tau_k)$ and in the latter $(\pi(z_k))$. Balanced sampling is easy to implement and often can reduce the variance substantially. In the recursive implementation for filtering, we can combine the resampling step at the end of the current iteration and the selection of the index at the beginning of the next iteration into a single selection of indices. In order to keep the notation simple, we discuss the ideas in the context of the resampling step only. We denote the number of times the index $k$ is selected by $N_k$. For random resampling with replacement, these multiplicities $(N_k)$ are multinomial $(N, (\pi(z_k)))$. Here, we look for more systematic sampling procedures. We require that $\sum N_k \equiv N$ and that sampling is unbiased, that is,

$$\mathbf{E}[N_j | z_1, \ldots, z_R] = N\pi(z_j).$$

Then the estimator

$$\hat{m}(\psi) = \frac{1}{N} \sum_{j=1}^{R} \psi(Z_j) N_j$$

has the same expected value as the usual importance sampling estimator

$$\tilde{m}(\psi) = \sum_{j=1}^{R} \psi(Z_j)\pi(Z_j).$$

Its variance can be written as

$$\mathrm{Var}(\hat{m}(\psi)) = \mathrm{Var}\left(\sum_{j=1}^{R} \psi(Z_j)\pi(Z_j)\right) + \frac{1}{N^2}\mathbf{E}\left[\sum_{i,j} \psi(Z_i)\psi(Z_j) C_R(i,j)\right],$$

where $C_R(i,j)$ is the conditional covariance of $N_i$ and $N_j$. The first term is the variance of the usual importance sampling estimator and the second term is the additional variability due to the resampling step. The advantage of resampling becomes apparent only when we consider several time steps: Without resampling, the recursive filter sample would quickly degenerate, that is, practically all the weights would be given to very few values. Resampling splits the particles with large weights into several independent ones and kills some of the particles with very small weights. Nevertheless, we should try to minimize the additional variability introduced by resampling. Because it is not known in advance which functions $\psi$ will be of interest, we consider the supremum over all (bounded) test functions $\psi$.

With multinomial $N_j$'s, we have

$$\begin{aligned}\sum_{i,j} \psi(z_i)\psi(z_j) C_R(i,j) &= N\left(\sum_i \psi(z_i)^2 \pi(z_i) - \left(\sum_i \psi(z_i)\pi(z_i)\right)^2\right) \\ &\le N \sup_x \psi(x)^2.\end{aligned}$$

Hence, resampling randomly with replacement can guarantee that the effect of resampling disappears asymptotically.

Several methods have been proposed which reduce the (conditional) variances $C_R(i,i)$. Residual sampling [22] takes

$$N_i = [N\pi(z_i)] + N_i', \qquad (N_i') \sim \text{multinomial}(N', (\pi'(z_i))),$$

where $[x]$ denotes the integer part of $x$ and

$$N' = N - \sum_i [N\pi(z_i)], \qquad \pi'(z_i) = \frac{N\pi(z_i) - [N\pi(z_i)]}{N'}.$$

This reduces $\sum_{i,j} \psi(z_i)\psi(z_j)C_R(i,j)$ by the factor $N'/N$. Intuitively, we expect the fractional part $N\pi(z_i) - [N\pi(z_i)]$ to be uniform on $(0,1)$, leading to an average reduction by a factor of two.

The variance $C_R(k,k)$ is minimal iff $N_k$ is equal to one of the two integers closest to $N\pi(z_k)$; see [6]. Together with the condition $\mathbf{E}[N_k] = N\pi(z_k)$, this determines the marginal distribution of $N_k$. Crisan and Lyons [6] show that then also the expected relative entropy between the empirical distribution $(N_k/N)$ and the target $(\pi(z_k))$ is minimal. There are at least two algorithms such that all $N_k$'s differ by less than one from $N\pi(z_k)$. The following one goes at least back to Whitley [28] and has been rediscovered by Carpenter, Clifford and Fearnhead [2]:

$$N_{j_k} = \left| \left[ N \sum_{i=1}^{k-1} \pi(z_{j_i}) + U, N \sum_{i=1}^{k} \pi(z_{j_i}) + U \right) \cap \{1, 2, \ldots, N\} \right|, \tag{12}$$

where $(j_1, j_2, \ldots, j_R)$ is a random permutation of $(1, 2, \ldots, R)$, $U$ is uniform on $[0,1)$, and the absolute value of a finite set denotes the number of elements in this set.

The second algorithm has been proposed by Crisan, Del Moral and Lyons [5]; see also [4]. One chooses an arbitrary binary tree with $R$ leaves, labelled as $1, \ldots, R$, and one propagates $N$ particles from the root in a specific way down the tree. The value of $N_j$ is then the number of particles ending at leaf $j$. In order to describe the propagation, we identify a node $\alpha$ with the subset of $\{1, \ldots, R\}$ that consists of the leaves connected to $\alpha$. Furthermore, we denote by $N_\alpha$ the number of particles that pass through a node $\alpha$. The expected value of $N_\alpha$ must then be equal to $\mu_\alpha = N \sum_{j \in \alpha} \pi_j$. The splitting at each node is done such that $N_\alpha$ differs by less than one from $\mu_\alpha$ and $\mathbf{E}[N_\alpha] = \mu_\alpha$. It is easy to see that this can be achieved: Each split is either deterministic or chooses between two possibilities with given probabilities. Decisions at different nodes are made independently.

However, by minimizing $C_R(k,k)$, we usually introduce strong dependence between different $N_j$'s, and the effects of this are hard to control. Trivially,

$|C_R(i,j)| \le 1/4$, but the bound

$$\sum_{i,j} \psi(z_i)\psi(z_j)C_R(i,j) \le \frac{N^2}{4} \sup_x \psi(x)^2$$

contains no useful information because it does not even allow one to conclude that the additional uncertainty due to resampling disappears asymptotically.

With Whitley's algorithm [28], $C_R(i,j)$ can be either positive or negative. Since we know nothing about the sign of $\psi(z_k)$, I do not see how one could obtain a better worst case bound. Still, the following lemma supports the conjecture that, on average, the algorithm (12) will behave well.

LEMMA 2. *For arbitrary probabilities* $(\pi_i)$ *and arbitrary* $N$, *consider the random variables*

$$N_j = \left| \left[ N \sum_{i=1}^{j-1} \pi_i + U, N \sum_{i=1}^{j} \pi_i + U \right) \cap \{1, 2, \dots, N\} \right|,$$

*where* $U$ *is uniform on* $(0,1)$ [*this is the algorithm* (12) *without the additional permutation*]. *Then for any* $j < k$, $\operatorname{Cov}(N_j, N_k)$ *depends only on* $r_l = N\pi_j \bmod 1$, $r_u = N\pi_k \bmod 1$ *and* $r_m = N \sum_{i=j+1}^{k-1} \pi_i \bmod 1$, *an explicit expression being given in the proof. Moreover, the average of this covariance with respect to the uniform distribution on* $(0,1)$ *for* $r_m$ *is zero for all values* $r_l$ *and* $r_u$.

PROOF. Because shifting a uniform random variable modulo 1 does not change the distribution, we may assume that $j = 1$. Moreover, it is clear that only the fractional parts $r_l, r_m, r_u$ matter. If we put $M_j = N_j - [N\pi_j]$ and $M_k = N_k - [N\pi_k]$, we obtain therefore

$$\begin{aligned} \mathbf{E}[M_j M_k] = {} & \mathbf{P}[U \in (0, r_l) \cap (r_l + r_m - 1, r_l + r_m + r_u - 1)] \\ & + \mathbf{P}[U \in (0, r_l) \cap (r_l + r_m - 2, r_l + r_m + r_u - 2)]. \end{aligned}$$

It is easy to evaluate the right-hand side by distinguishing different cases:

$$\mathbf{E}[M_j M_k] = \begin{cases} (r_l + r_m + r_u - 1)^+, & (r_l + r_m \le 1, r_m + r_u \le 1), \\ r_u, & (r_l + r_m > 1, r_m + r_u \le 1), \\ r_l, & (r_l + r_m \le 1, r_m + r_u > 1), \\ 1 - r_m, & (r_l + r_m > 1, r_m + r_u > 1, \\ & \qquad r_l + r_m + r_u \le 2), \\ r_l + r_u - 1, & (r_l + r_m + r_u > 2). \end{cases}$$

It is also easy to show that, by integrating over $r_m \in (0,1)$, we obtain $r_l r_u$ in all cases. $\square$

Although the method is unbiased and has minimal variance without randomizing the order of the values, it seems wise to do so since it is computationally easy and we expect it to make the values $r_m$ approximately uniform.

With the algorithm of Crisan, Del Moral and Lyons [5] we have control over the sign of $C_R(j,k)$.

LEMMA 3. *For the tree-based algorithm, we have, for arbitrary nondecreasing functions* $f_1, \dots, f_R$,

$$\mathbf{E}\left[\prod_{j=1}^{R} f_j(N_j)\right] \le \prod_{j=1}^{R} \mathbf{E}[f_j(N_j)].$$

*In particular, the covariances* $C_R(i,j)$ *are negative for* $i \ne j$.

PROOF. Denote the two nodes connected directly to the root by $\alpha$ and $\beta$. Because the particles are propagated independently in the two subtrees,

$$\mathbf{E}\left[\prod_{j=1}^{R} f_j(N_j) | N_\alpha, N_\beta\right] = \mathbf{E}\left[\prod_{j\in\alpha} f_j(N_j) | N_\alpha\right] \mathbf{E}\left[\prod_{j\in\beta} f_j(N_j) | N_\beta\right].$$

Furthermore,

$$\mathbf{E}\left[\prod_{j\in\alpha} f_j(N_j) | N_\alpha = [\mu_\alpha]\right] \le \mathbf{E}\left[\prod_{j\in\alpha} f_j(N_j) | N_\alpha = [\mu_\alpha] + 1\right]$$

since we can propagate first $[\mu_\alpha]$ particles and afterward an additional particle. Because $N_\alpha$ and $N_\beta$ are negatively dependent, we obtain

$$\mathbf{E}\left[\prod_{j=1}^{R} f_j(N_j)\right] \le \mathbf{E}\left[\prod_{j\in\alpha} f_j(N_j)\right] \mathbf{E}\left[\prod_{j\in\beta} f_j(N_j)\right].$$

The proof proceeds now recursively, by considering in the next step each factor separately and conditioning on the number of particles one level lower. □

This lemma implies that the additional variance due to resampling is reduced by a factor of at least two compared to multinomial sampling:

LEMMA 4. *For the tree-based algorithm described,*

$$\sum_{i,j} \psi(z_i)\psi(z_j) C_R(i,j) \le \frac{1}{2} \sum_i \psi(z_i)^2 \le \frac{N}{2} \sup_x \psi(x)^2.$$

PROOF. Write $\psi$ as the difference of positive and negative parts and use Cauchy–Schwarz for the covariance between positive and negative parts; see also [4], page 31. □

For later use we formulate and prove the following large deviation inequality:

LEMMA 5. *For the tree-based algorithm,*

$$\sup_{A\subset\{1,\dots,R\}} \mathbf{P}\left[\left|\sum_{j\in A}(N_j - N\pi_j)\right| \ge \varepsilon\right] \le 2\exp(-4\varepsilon^2/R).$$

PROOF. For any $t > 0$, we have

$$\mathbf{P}\left[\sum_{j\in A}(N_j - N\pi_j) \ge \varepsilon\right] \le \exp(-t\varepsilon)\mathbf{E}\left[\exp\left(t\sum_{j\in A}(N_j - N\pi_j)\right)\right].$$

By Lemma 3,

$$\mathbf{E}\left[\exp\left(t\sum_{j\in A}(N_j - N\pi_j)\right)\right] \le \prod_{j\in A}\mathbf{E}[\exp(t(N_j - N\pi_j))].$$

Because $N_j$ takes only two values,

$$\mathbf{E}[\exp(t(N_j - N\pi_j))] = \exp(-tr_j)(1 - r_j + \exp(t)r_j),$$

where $r_j = N\pi_j - [N\pi_j]$. By the standard argument in the proof of Hoeffding's inequality, the right-hand side can be bounded by $\exp(t^2/8)$. Hence,

$$\mathbf{P}\left[\sum_{j\in A}(N_j - N\pi_j) \ge \varepsilon\right] \le \exp(-t\varepsilon + |A|t^2/8),$$

which is minimal for $t = 4\varepsilon/|A|$. The probability of a deviation less than or equal to $-\varepsilon$ can be bounded by the same expression. The lemma then follows because we may assume $|A| \le R/2$ (the deviations for $A$ and $A^c$ differ only in sign). □

3.3. *Accept–reject versus sampling importance resampling.* Generally speaking, the computational effort for rejection sampling is greater than for sampling importance resampling. By how much depends, however, on the specific situation. Note that with the auxiliary variable idea of Pitt and Shephard [23], it is possible to use the rejection method even in cases where the likelihood $b$ is unbounded, for example, in the stochastic volatility model of Section 3.1.1 with $y = 0$. For both methods one needs to find proposal densities $\rho(j, x)$ that approximate $a(j, x)b(x)$, but for rejection sampling one needs, in addition, an upper bound for $a(j,x)b(x)/\rho(j,x)$, which can be difficult in high dimensions.

Usually, a large empirical variance of the inclusion probabilities $\pi(z_k)$ is taken as an indication that the error of sampling importance resampling is large. However, a low variance does not guarantee a low error. When the true filter density is bimodal and if the proposal represents only one mode well, then the inclusion probabilities are fairly balanced unless the sample size is huge. If we are unable to compute the modes of the filter density, then rejection sampling is presumably the only way to obtain some guarantee for the algorithm in such a case.

The results of the next section allow some theoretical comparison of rejection and sampling importance resampling methods. We will show in Section 3.1.1 below that rejection sampling has a smaller asymptotic variance than the standard sampling importance resampling algorithm. Another relevant question is whether the errors of the methods depend on $\sup_x b_t(x, y_t)$ (if they do, then it is not clear how much one gains by an algorithm which does not need a bound on this supremum). For the rejection method, both the exponential bounds in the law of large numbers and the asymptotic variance do not depend on $b_t$ at all as long as the condition (20) is satisfied. For sampling importance resampling, our best bound in the law of large numbers depends on $\sup_x b_t(x, y_t)$ because of Lemma 9. The bound on the asymptotic variance does not involve the supremum of $b_t$, but a certain $L_2$-norm of $b_t$.

3.4. *Computation of the likelihood.* Combining (3) and (6), we see that

$$p(y_t|y_{1:t-1}) \approx \sum_{j=1}^{N} \int \frac{1}{N} a_t(x_{j,t-1}, x) b_t(x, y_t)\, d\mu(x),$$

which is in the short notation of this section equal to $\sum \beta_j/N$. If we use $\tau_j \rho(j, x)$ as our proposal, then the usual importance sampling estimator of $p(y_t|y_{1:t-1})$ is

$$\widehat{p}(y_t|y_{1:t-1}) = \frac{1}{NR} \sum_{k=1}^{R} \frac{b(z_k)a(j_k, z_k)}{\tau_{j_k}\rho(j_k, z_k)}.$$

3.5. *Monte Carlo backward smoothing.* There is a similar recursive simulation method that generates samples from the conditional distribution of $X_{0:T}$ given $Y_{1:T} = y_{1:T}$. At time $T$, we use the recursive filter sample $x_{j,T}^{sm} = x_{j,T}$. We then proceed backward in time, using (4) together with an approximation of $f_{t|t}$. In order to avoid problems with discreteness, we recommend use of (6) as in [18], instead of replacing $f_{t|t}$ by the empirical distribution of the particles at time $t$ as in [14]. This means that we generate $x_{j,t}^{sm}$ from $x_{j,t+1}^{sm}$ and $(x_{i,t-1})$ by simulating from the density proportional to

$$a_{t+1}(x, x_{j,t+1}^{sm}) b_t(x, y_t) \frac{1}{N} \sum_{i=1}^{N} a_t(x_{i,t-1}, x). \tag{13}$$

[At time $t = 0$, we will use the density proportional to $a_1(x, x_{j,1}^{sm}) a_0(x)$.] Clearly this has the same structure as (7) and so the same methods as discussed before apply in principle. However, we need one value from the density (13) for each $j$ and thus sampling importance resampling does not seem to be useful here. For the same reason, care is needed when using the mixture index as an auxiliary variable. Since sampling from $(\tau_i)$ typically involves computing the partial sums of the $\tau_i$'s, one should use the same distribution $(\tau_i)$ for all $j$. Then the computational cost of the approach is $O(TN)$ and thus at least comparable to a standard MCMC method. The main disadvantage of this approach is that we have to store all the filter samples.

**4. Theoretical properties.** In this section we analyze the convergence of the approximation $f_{t|t}^N$ to the true filtering density $f_{t|t}$. We will hold the observations $y_{1:t}$ fixed and drop them from the notation. In particular, we do not make any assumption about how the observations were obtained. The true filtering densities $f_{t|t}$ are then deterministic, but the approximations $f_{t|t}^N$ are still random since their computation involves random sampling. All expectations and probabilities in this section concern the randomness of the Monte Carlo methods, and not the randomness of the state space model. We assume throughout that $X_t$ takes its values in a complete, separable metric space equiped with the Borel $\sigma$-field, and we denote the metric on this state space by $d(\cdot,\cdot)$.

The operator notation for recursive Monte Carlo filters introduced in Section 2.3 will be used extensively. In addition, we denote by $E_N(f)$ the empirical distribution of a sample of size $N$ from $f$. Then the approximate filter density is

$$f_{t|t}^N = B(A_t^* E_N(f_{t-1|t-1}^N), b_t(\cdot, y_t)),$$

and by (6) and (5) it has to be compared with

$$f_{t|t} = B(A_t^* f_{t-1|t-1}, b_t(\cdot, y_t)).$$

In the first two sections we present two approaches for showing convergence of $f_{t|t}^N$ to $f_{t|t}$ as $N \to \infty$. We measure the error by the $L_1$-distance between densities, see, for example, [9], Chapter 1, which can be written in several equivalent forms:

$$\begin{aligned}\|f-g\|_1 &= \int |f(x)-g(x)|\,d\mu(x) = 2\int (f(x)-g(x))^+\,d\mu(x) \\ &= 2\sup_C |P_f[C]-P_g[C]| = 2\int (f(x)-\min(f(x),g(x)))\,d\mu(x)\end{aligned} \tag{14}$$

($x^+$ denotes the positive part of $x$). Clearly, if $\|f^N_{t|t} - f_{t|t}\|_1$ converges to zero in probability or almost surely, then for any bounded function $\psi$ on the state space, the law of large number holds:

$$\frac{1}{N}\sum_{j=1}^N \psi(x_{j,t}) \longrightarrow \int \psi(x) f_{t|t}(x)\,d\mu(x)$$

in probability or almost surely. In the third section we show the corresponding central limit theorem.

4.1. *Stepwise error propagation.* The obvious first attempt to show convergence uses the decomposition

$$\begin{aligned}f^N_{t|t} - f_{t|t} &= B(A^*_t E_N(f^N_{t-1|t-1}), b_t) - B(A^*_t f^N_{t-1|t-1}, b_t) \\ &\quad + B(A^*_t f^N_{t-1|t-1}, b_t) - B(A^*_t f_{t-1|t-1}, b_t).\end{aligned} \tag{15}$$

The first term is the error due to sampling at time $t-1$ (propagated once) and the second term is the propagation of the error at time $t-1$. For a recursive inequality for $\|f^N_{t|t} - f_{t|t}\|_1$, we have to study the Lipschitz-continuity of Bayes and Markov operators with respect to the $L_1$-distance and to control the sampling error.

The continuity of Markov operators is well known; see [10], Section 3.

LEMMA 6. *We have*

$$\|A^*f - A^*g\|_1 \le \rho(A^*)\|f-g\|_1,$$

*where*

$$\rho(A^*) = \tfrac{1}{2}\sup_{x,x'} \|a(x,\cdot) - a(x',\cdot)\|_1 \le 1.$$

Note that, for a compact state space, the Markov operator is typically contracting.

The continuity of Bayes' formula with respect to the prior is more problematic. We have, see [20], Lemma 3.6(i), the following:

LEMMA 7.

$$\|B(f,b) - B(g,b)\|_1 \le \beta(f,b)\|f-g\|_1,$$

*where*

$$\beta(f,b) = \frac{\sup_x b(x)}{\int b(x)f(x)\,d\mu(x)} \in \left(1, \frac{\sup_x b(x)}{\inf_x b(x)}\right].$$

The difficulty is that this bound cannot be improved in general. Lemma 3.6(ii) from [20] shows that the Bayes operator is not contracting for any $f$, at least for some "directions" $g$.

Finally, we have the following bound on sampling errors:

LEMMA 8. *If $x \to a(x,\cdot)$ is continuous with respect to the $L_1$-norm, then under i.i.d. sampling from $g$,*

$$\mathbf{P}[\|A^* E_N(g) - A^* g\|_1 > \varepsilon] \overset{N\to\infty}{\longrightarrow} 0$$

*exponentially fast in $N$ for any $\varepsilon > 0$. The convergence is uniform for all $g$ such that $\int_K g\,d\mu \ge 1 - \varepsilon/6$ for some fixed compact set $K$.*

PROOF. The proof follows closely the arguments in [9], Chapter 3. We denote by $\mu_N$ the empirical distribution $E_N(g)$ and by $\mu_g$ the distribution $g(x)\,d\mu(x)$.

Let $\varepsilon > 0$ be given. Choose a compact $K$ such that $\mu_g(K) \ge 1 - \varepsilon/6$. Next, choose $\delta$ such that $\|a(x,\cdot) - a(x',\cdot)\|_1 \le \varepsilon/6$ for all $x, x' \in K$ with $d(x,x') \le \delta$. Then choose a partition $\{B_1, \dots, B_J\}$ of $K$ such that each $B_j$ has diameter at most $\delta$ and choose a point $z_j$ in $B_j$ for each $j$. Finally, put $B_0 = K^c$. Then

$$\begin{aligned}
&\left\| \frac{1}{N}\sum_{i=1}^N a(x_i,\cdot) - \sum_{j=1}^J \mu_N(B_j) a(z_j,\cdot) \right\|_1 \\
&\qquad = \int \left| \frac{1}{N}\sum_{i=1}^N \mathbb{1}_{B_0}(x_i) a(x_i,x) + \sum_{j=1}^J \frac{1}{N}\sum_{i=1}^N \mathbb{1}_{B_j}(x_i)(a(x_i,x) - a(z_j,x)) \right| d\mu(x) \\
&\qquad \le \mu_N(B_0) + \sum_{j=1}^J \frac{1}{N}\sum_{i=1}^N \mathbb{1}_{B_j}(x_i) \int |a(x_i,x) - a(z_j,x)|\,d\mu(x) \\
&\qquad \le |\mu_N(B_0) - \mu_g(B_0)| + \frac{\varepsilon}{3}.
\end{aligned}$$

Similarly, we obtain

$$\left\| \int a(x,\cdot) g(x)\,d\mu(x) - \sum_{j=1}^J \mu_g(B_j) a(z_j,\cdot) \right\|_1 \le \frac{\varepsilon}{3}.$$

Finally,

$$\left\|\sum_{j=1}^{J}\mu_g(B_j)a(z_j,\cdot)-\sum_{j=1}^{J}\mu_N(B_j)a(z_j,\cdot)\right\|_1\le\sum_{j=1}^{J}|\mu_N(B_j)-\mu_g(B_j)|.$$

Taking these three inequalities together, we obtain

$$\|A^*E_N(g)-A^*g\|_1\le\frac{2\varepsilon}{3}+\sum_{j=0}^{J}|\mu_N(B_j)-\mu_g(B_j)|.$$

Hence, the large deviation estimate for the multinomial distribution,

$$\mathbf{P}\left[\sum_{j=0}^{J}|\mu_N(B_j)-\mu_g(B_j)|>\frac{\varepsilon}{3}\right]\le 2^{J+2}\exp(-N\varepsilon^2/18)$$

([9], Theorem 3.2), implies

$$\mathbf{P}[\|A^*E_N(g)-A^*g\|_1>\varepsilon]\le 2^{J+2}\exp(-N\varepsilon^2/18).$$

From this, the lemma follows (note that, once $K$ is fixed, $J$ depends only on the transition kernel $a$ and not on $g$). □

THEOREM 1. *If $x\to a_t(x,\cdot)$ is continuous and if for all $t$, all $x$ and all $y$,*

$$0<b_t(x,y)\le C(t,y)<\infty,$$

*then for all $t$ and all $y_{1:t}$,*

$$\|f^N_{t|t}-f_{t|t}\|_1\longrightarrow 0$$

*in probability as $N\to\infty$.*

PROOF. The proof proceeds by induction on $t$. For $t=0$, there is nothing to prove because $f^N_{0|0}=f_{0|0}=a_0$. From Lemmas 6 and 7, it follows that

$$\begin{aligned}&\|B(A_t^*f^N_{t-1|t-1},b_t)-B(A_t^*f_{t-1|t-1},b_t)\|_1\\&\quad\le\frac{C(t,y_t)}{p(y_t|y_{1:t-1})}\|f^N_{t-1|t-1}-f_{t-1|t-1}\|_1\le\varepsilon\end{aligned}$$

if

$$\|f^N_{t-1|t-1}-f_{t-1|t-1}\|_1\le\varepsilon\frac{p(y_t|y_{1:t-1})}{C(t,y_t)}=:\delta. \tag{16}$$

By the induction assumption, there is an $N_1$ such that, for $N>N_1$, (16) holds with probability at least $1-\eta$.

In order to bound the first term in (15), some care is needed when applying the bounds provided by Lemmas 7 and 8 with $f^N_{t-1|t-1}$, which is random. We have to show that when (16) holds, we can obtain bounds which depend only on $f_{t-1|t-1}$. Note first that

$$\int b_t(x,y_t)(A_t^* f^N_{t-1|t-1}(x) - A_t^* f_{t-1|t-1}(x))\,d\mu(x)$$
$$\geq -\tfrac{1}{2}C(t,y_t)\|f^N_{t-1|t-1} - f_{t-1|t-1}\|_1.$$

Hence, if (16) is satisfied,

$$\int b_t(x,y_t)A_t^* f^N_{t-1|t-1}(x)\,d\mu(x) \geq (1-\varepsilon/2)p(y_t|y_{1:t-1}) \geq \tfrac{1}{2}p(y_t|y_{1:t-1})$$

and, therefore, by Lemma 7, also

$$\|B(A_t^* E_N(f^N_{t-1|t-1}), b_t) - B(A_t^* f^N_{t-1|t-1}, b_t)\|_1$$
$$\leq \frac{2C(t,y_t)}{p(y_t|y_{1:t-1})}\|A_t^* E_N(f^N_{t-1|t-1}) - A_t^* f^N_{t-1|t-1}\|_1.$$

Next we observe that, if $K$ is compact such that $\int_K f_{t-1|t-1}\,d\mu \geq 1-\delta/2$ and if (16) holds, then $\int_K f^N_{t-1|t-1}\,d\mu \geq 1-\delta$. Therefore, by Lemma 8 we can find $N_2$ such that, for $N > N_2$,

$$\|A_t^* E_N(f^N_{t-1|t-1}) - A_t^* f^N_{t-1|t-1}\|_1 \leq 6\delta \tag{17}$$

holds with probability at least $1-\eta$. Collecting all the bounds shows that, for $N > \max(N_1, N_2)$,

$$\|f^N_{t|t} - f_{t|t}\|_1 \leq 13\varepsilon$$

with probability at least $1-2\eta$. □

The conditions of this theorem are weak. However, the arguments in the proof require $\|f^N_{t-1|t-1} - f_{t-1|t-1}\|_1$ to be smaller than $\|f^N_{t|t} - f_{t|t}\|_1$. This means that the required sample size $N$ grows with $t$. It is easy to see that, in general, $N$ has to grow exponentially with $t$, and, thus, from a practical point of view, the theorem is not of great use. Strengthening the assumptions by, for instance, assuming a compact state space, does not help because by Lemma 3.6(ii) from [20], the Bayes operator is expanding. Hence, for a more useful result, we need a different approach which is provided in the next section.

4.1.1. *Sampling errors for sampling importance resampling.* The results so far have assumed that the Monte Carlo filter uses i.i.d. samples of $f_{t|t}^N$, which means using the accept–reject method (with or without auxiliary variables). It does not cover sampling importance resampling. In order to extend the results above, we need to adapt Lemma 8 to the different sampling method.

LEMMA 9. *Let $g$ have the form $g = B(h,b)$ and let $(x_i, N_i)$ be a sampling importance resample from $g$, that is, $(x_i)$ is an i.i.d. sample from $h$ and the $N_i$'s are the multiplicities in the resampling step which uses probabilities $\pi_i \propto b(x_i)$. Assume that $x \to a(x,\cdot)$ is continuous for the $L_1$-norm, that $\sup b(x)/\int b(x)h(x)\,d\mu(x) < \infty$ and that*

$$\sup_{J\subset\{1,2,\ldots,N\}} \mathbf{P}\left[\left|\sum_{j\in J}\left(\frac{N_j}{N} - \pi_j\right)\right| > \varepsilon\right] \le c_1 \exp(-c_2 N\varepsilon^2).$$

*Then*

$$\mathbf{P}[\|A^* E_N(g) - A^* g\|_1 > \varepsilon] \overset{N\to\infty}{\longrightarrow} 0$$

*exponentially fast in $N$ for any $\varepsilon > 0$.*

PROOF. The assumption of i.i.d. sampling was used in the proof of Lemma 8 only to obtain an exponential bound for

$$\mathbf{P}\left[\sum_{j=0}^{J} |\mu_N(B_j) - \mu_g(B_j)| > \frac{\varepsilon}{3}\right].$$

Hence, we have to obtain such a bound by different arguments. By Scheffé's theorem and Bonferroni's inequality, we have

$$\mathbf{P}\left[\sum_{j=0}^{J} |\mu_N(B_j) - \mu_g(B_j)| > \frac{\varepsilon}{3}\right] \le 2^{J+1} \sup_B \mathbf{P}\left[|\mu_N(B) - \mu_g(B)| > \frac{\varepsilon}{6}\right],$$

where the supremum is taken over all sets $B$ in the $\sigma$-field generated by $(B_0, B_1, \ldots, B_J)$. We can decompose via

$$\begin{aligned}
&\mu_N(B) - \mu_g(B)\\
&\quad = \sum_{i=1}^{N}\left(\frac{N_i}{N} - \pi_i\right)\mathbb{1}_B(x_i)\\
&\qquad - \frac{1}{\int b(x)h(x)\,d\mu(x)}\left(\frac{1}{N}\sum_{i=1}^{N} b(x_i) - \int b(x)h(x)\,d\mu(x)\right)\frac{\sum_{i=1}^{N} b(x_i)\mathbb{1}_B(x_i)}{\sum_{i=1}^{N} b(x_i)}\\
&\qquad + \frac{1}{\int b(x)h(x)\,d\mu(x)}\left(\frac{1}{N}\sum_{i=1}^{N} b(x_i)\mathbb{1}_B(x_i) - \int_B b(x)h(x)\,d\mu(x)\right).
\end{aligned}$$

The assumption on the resampling method gives an exponential bound for the probability that the first term is larger than $\varepsilon/18$. Hoeffding's inequality provides analogous bounds for the second and third terms. □

Applying this lemma with $h = N^{-1}\sum_j a(x_{j,t-2},\cdot)$ and $b = b_{t-1}(\cdot, y_{t-1})$, we obtain the analogue of Theorem 1. The arguments in the proof of this theorem show that in this case $b(x)/\int b(x)h(x)\,d\mu(x)$ is bounded.

4.2. *Analysis based on considering several steps.* Clearly, we can look at error propagation over more than one time step. If we define

$$K_{s,t}(f) = K_{s+1,t}(B(A^*_{s+1}f, b_{s+1})) \qquad (s < t),\ \ K_{t,t}(f) = f,$$

then, for any $s < t$, $f_{t|t} = K_{s,t}(f_{s|s})$ and, hence,

$$\begin{aligned} f^N_{t|t} - f_{t|t} = \sum_{r=s+1}^{t} (K_{r,t}(B(A^*_r E_N(f^N_{r-1|r-1}), b_r)) \\ - K_{r,t}(B(A^*_r f^N_{r-1|r-1}, b_r))) \\ + K_{s,t}(f^N_{s|s}) - K_{s,t}(f_{s|s}). \end{aligned} \tag{18}$$

Here the last difference is the error at time $s$ propagated over $t-s$ steps. The other differences are the errors due to sampling at time $r-1$, propagated over $t-r+1$ steps.

This is only useful if we can give a bound on the error propagated over $k$ steps which is better than the sum over $k$ single steps. It is possible because an alternative way to get from $f_{s|s}$ to $f_{t|t}$ is to apply first the Bayes operator once with likelihood equal to the conditional density of $y_{s+1:t}$ given $x_s$, followed by $t-s$ Markov operators for the conditional transitions from $x_r$ to $x_{r+1}$ given $y_{r+1:t}$. The contractivity of the Markov operators can then beat the expansion of the Bayes operator. It requires, however, a uniform nontrivial upper bound for the contraction coefficient of the conditional chain given $y_{r+1:t}$, and for this, we need the following condition: There are probability densities $h_t$ and two constants $0 < c_a < C_a < \infty$ such that, for all $x$ and $x'$,

$$c_a\,h_t(x) \le a_t(x', x) \le C_a\,h_t(x). \tag{19}$$

Condition (19) on $a_t$ is reasonable when the state space is compact, although it is slightly stronger than uniform ergodicity. Using (14), we see that the lower bound of (19) alone implies $\rho(A^*_t) \le 1 - c_a$ and thus also uniform ergodicity. Condition (19) includes even some examples with unbounded state space. For instance, (19) holds for the model

$$X_t = g(X_{t-1}) + V_t$$

if $g$ is bounded and $V_t$ has a density whose logarithm is uniformly Lipschitz continuous. This is satisfied for most heavy-tailed distributions, but not for the Gaussian. For Gaussian $V_t$, (19) is false: There is no density $h_t$ such that the two bounds in (19) hold simultaneously. We thus have an example of a uniformly ergodic chain that we cannot treat with our arguments.

Concerning $b_t$, there is an almost minimal condition, namely,

$$0 < \int b_t(x, y_t) h_t(x)\, d\mu(x) < \infty \tag{20}$$

for all $t$ and all $y_t$. Some arguments become much simpler, however, if we replace (20) by

$$C_b := \sup_{t,x,x',y} \frac{b_t(x, y)}{b_t(x', y)} < \infty. \tag{21}$$

The following lemma shows that, under condition (19), the error propagated over several steps decreases exponentially. Many versions of this exponential forgetting of the initial conditions of the filter have appeared in the literature; see, for example, [7, 8, 21] and the references given there. We use the version of [11], Lemma 1.

LEMMA 10. *Assume conditions* (19) *and* (20). *Then for any two densities* $f$ *and* $g$ *and any* $s < t$ *we have*

$$\|K_{s,t}(f) - K_{s,t}(g)\|_1 \le \frac{1}{\gamma_a}(1 - \gamma_a)^{t-s}\|f - g\|_1,$$

*where* $\gamma_a = c_a/C_a$.

PROOF. As already mentioned, we write $K_{s,t}$ as the composition of one Bayes operator and $t - s$ Markov operators. The likelihood in the Bayes operator is equal to the conditional density of $y_{s+1:t}$ given $x_s$. It satisfies the recursion

$$p(y_{s+1:t}|x_s) = \int a_{s+1}(x_s, x_{s+1}) b_{s+1}(x_{s+1}, y_{s+1}) p(y_{s+2:t}|x_{s+1})\, d\mu(x_{s+1}).$$

Hence, by conditions (19) and (20) and an induction argument,

$$\frac{\sup_{x_s} p(y_{s+1:t}|x_s)}{\inf_{x_s} p(y_{s+1:t}|x_s)} \le \frac{1}{\gamma_a}, \tag{22}$$

which is, by Lemma 7, the maximal expansion by the Bayes operator. The Markov operators have transition densities

$$p(x_r|x_{r-1}, y_{r:t}) = \frac{a_r(x_{r-1}, x_r) b_r(x_r, y_r) p(y_{r+1:t}|x_r)}{p(y_{r:t}|x_{r-1})},$$

which are bounded below by

$$\gamma_a \frac{h_r(x_r) b_r(x_r, y_r) p(y_{r+1:t}|x_r)}{\int h_r(x_r) b_r(x_r, y_r) p(y_{r+1:t}|x_r)\, d\mu(x_r)}.$$

The right-hand side is $\gamma_a$ times a density that does not depend on $x_{r-1}$. Hence, by Lemma 6 and (14), each Markov operator contracts at least by $(1-\gamma_a)$. □

THEOREM 2. *Assume that the transition densities $a_t$ are the same for all $t$, that they are continuous in the $L_1$-norm and satisfy (19), and that (21) holds. Then to any $\varepsilon > 0$, there are constants $c_1$ and $c_2$ such that, for all $t$ and all $N$,*

$$\mathbf{P}[\|f^N_{t|t} - f_{t|t}\|_1 > \varepsilon] \le c_1 \exp(-c_2 N).$$

PROOF. Because $a_t$ and thus also $A^*_t$ are the same for all $t$, we drop the time index during this proof. Let $\varepsilon > 0$ be given. Choose $k$ such that

$$\frac{2}{\gamma_a}(1-\gamma_a)^k \le \varepsilon.$$

Assume first that $k < t$. Because the $L_1$-distance between densities is at most 2, we obtain, in this case from the decomposition (18) with $s = t - k$ and Lemmas 10 and 7,

$$\begin{aligned}
&\|f^N_{t|t} - f_{t|t}\|_1 \\
&\qquad \le \frac{1}{\gamma_a} \sum_{r=t-k+1}^{t} (1-\gamma_a)^{t-r} \|B(A^* E_N(f^N_{r-1|r-1}), b_r) \\
&\qquad\qquad\qquad\qquad - B(A^*(f^N_{r-1|r-1}), b_r)\|_1 + \varepsilon \\
&\qquad \le \frac{C_b}{\gamma_a} \sum_{r=t-k+1}^{t} (1-\gamma_a)^{t-r} \|A^* E_N(f^N_{r-1|r-1}) - A^* f^N_{r|r}\|_1 + \varepsilon.
\end{aligned}$$

If $k > t$, we obtain a similar result by considering the decomposition (18) with $s = 0$. (Because $f^N_{0|0} = f_{0|0} = a_0$, the $\varepsilon$ at the end is then absent.) Hence, if

$$\sup_{t-k \le r < t} \|A^* E_N(f^N_{r|r}) - A^* f^N_{r|r}\|_1 \le \varepsilon \tag{23}$$

holds, then, by the formula for a geometric series,

$$\|f^N_{t|t} - f_{t|t}\|_1 \le (C_b \gamma_a^{-2} + 1)\varepsilon.$$

We are now going to bound the probability that (23) occurs. Note that $\varepsilon$ and thus also $k$ are fixed. Because of Lemma 8, all we need to show is that the set of distributions $(f_{r|r}^N\,d\mu)$ is tight. By the definition of $f_{r|r}^N$ and by the conditions (19) and (21), we have

$$f_{r|r}^N(x) = \frac{\sum_{j=1}^N a(x_{j,r-1},x)b_r(x,y_r)}{\sum_{j=1}^N \int a(x_{j,r-1},x)b_r(x,y_r)\,d\mu(x)} \le C_b C_a h(x).$$

Clearly this implies the desired tightness. □

The important feature of the above theorem is that the same $N$ works for all times $t$. By Bonferroni's inequality, we obtain

$$\mathbf{P}\left[\sup_{t\le T}\|f_{t|t}^N - f_{t|t}\|_1 > \varepsilon\right] \le Tc_1\exp(-c_2 N).$$

Hence, it is sufficient to let $N$ increase logarithmically with the length of the series to guarantee uniform convergence of the filter approximation at all time points. It is not difficult to extend the above theorem to cases where the state transitions depend on $t$ as long as the continuity is uniform in $t$.

Condition (21) is used in the proof for bounding

$$\|B(A^*E_N(f_{r|r}^N),b_{r+1}) - B(A^*f_{r|r}^N,b_{r+1})\|_1$$

by applying Lemmas 7 and 8. The following lemma provides a direct way to bound the above distance by imposing only conditions on $a$, but assuming a compact state space.

Lemma 11. *Let $a$ be a transition density on a compact state space that satisfies (19) and*

$$\Delta(x',x) := \sup_{x''}\frac{|a(x,x'')-a(x',x'')|}{h(x'')} \to 0 \qquad [d(x,x')\to 0] \tag{24}$$

*with the same density $h$ as in (19). Then under i.i.d. sampling from $g$,*

$$\mathbf{P}[\|B(A^*E_N(g),b) - B(A^*g,b)\|_1 > \varepsilon] \overset{N\to\infty}{\longrightarrow} 0$$

*exponentially fast in $N$ for any $\varepsilon>0$, uniformly over all densities $g$ and all likelihoods $b$ with $0<\int h(x')b(x')\,d\mu(x')<\infty$.*

Proof. To make the notation more compact, we introduce

$$q(x',x) = \frac{a(x',x)b(x)}{\beta(x')}, \qquad \beta(x) = \int a(x,x')b(x')\,d\mu(x').$$

Then $q(x',x)$ is again a transition density and we can write

$$B(A^*E_N(g),b)(x) = \sum_{i=1}^N \frac{\beta(x_i)}{\sum_{k=1}^N \beta(x_k)} q(x_i,x)$$

and

$$B(A^*g,b)(x)=\int\frac{g(x')\beta(x')}{\int g(x'')\beta(x'')\,d\mu(x'')}q(x',x)\,d\mu(x').$$

The difference between these two expressions can thus be decomposed as

$$-\frac{1}{\int g(x)\beta(x)\,d\mu(x)}\left(\frac{1}{N}\sum_{j=1}^{N}\beta(x_j)-\int g(x)\beta(x)\,d\mu(x)\right)\frac{\sum_{j=1}^{N}\beta(x_j)q(x_j,x)}{\sum_{j=1}^{N}\beta(x_j)}$$
$$+\frac{1}{\int g(x)\beta(x)\,d\mu(x)}\left(\frac{1}{N}\sum_{j=1}^{N}\beta(x_j)q(x_j,x)\right.$$
$$\left.-\int g(x')\beta(x')q(x',x)\,d\mu(x')\right).$$

By assumption (19), we have

$$\gamma_a\le\frac{\beta(x)}{\int g(x)\beta(x)\,d\mu(x)}\le\gamma_a^{-1}.$$

Hence, it follows by Hoeffding's inequality that

$$\mathbf{P}\left[\left|N^{-1}\sum_{j=1}^{N}\frac{\beta(x_j)}{\int g(x)\beta(x)\,d\mu(x)}-1\right|>\varepsilon\right]\le 2\exp(-2N\varepsilon^2\gamma_a^2/(1-\gamma_a^2)^2).$$

Because the $L_1$-norm of $\sum_j\beta(x_j)q(x_j,x)/\sum_j\beta(x_j)$ is one, we have the same bound for the probability that the $L_1$-norm of the first term is greater than $\varepsilon$.

Assumption (24) allows us to control the continuity of $x\to\beta(x)q(x,\cdot)=a(x,\cdot)b(\cdot)$ with respect to the $L_1$-norm:

$$\|a(x,\cdot)b(\cdot)-a(x',\cdot)b(\cdot)\|_1\le\Delta(x,x')\int h(x'')b(x'')\,d\mu(x'').$$

Hence, the same argument as in Lemma 8 can be used to prove an exponential bound for the probability that the $L_1$-norm of the second term is greater than $\varepsilon$. □

By looking at the proof of Theorem 2, this lemma implies immediately the following:

THEOREM 3. *The claim of Theorem 2 is valid if the state space is compact, the transition densities do not depend on $t$ and (19), (20) and (24) hold.*

4.3. *Central limit theorems.* The goal of this section is to show by a simple induction argument that

$$\sqrt{N}\left(\frac{1}{N}\sum_{j=1}^{N}\psi_s(x_{j,s}) - \int \psi_s(x) f_{s|s}(x)\,d\mu(x)\right)_{0\le s\le t}$$

is asymptotically centered normal for any fixed $t$, any $y_{1:t}$ and functions $\psi_s$, $0 \le s \le t$, which are square integrable w.r.t. $f_{s|s}$. Del Moral and Miclo ([8], Corollary 20) have obtained a similar result, but we do not assume the $\psi_s$'s to be bounded nor the likelihood $b_t(\cdot, y_t)$ to be bounded away from zero.

Our argument proceeds by induction on the number $t$ of time steps. For $t = 0$, the result is obvious because $(x_{j,0})$ is an i.i.d. sample from $f_{0|0} = a_0$. The key idea for the induction step is to condition on $(x_{j,t-1})$. We first explain the argument heuristically. Introducing the notation

$$M_{N,t}(\psi) = \frac{1}{N}\sum_{j=1}^{N}\psi(x_{j,t}),$$

$$m_{N,t}(\psi) = \int \psi(x) f^N_{t|t}(x)\,d\mu(x),$$

$$m_t(\psi) = \int \psi(x) f_{t|t}(x)\,d\mu(x),$$

we can split

$$\begin{aligned} &\sqrt{N}(M_{N,t}(\psi) - m_t(\psi)) \\ &\quad = \sqrt{N}(M_{N,t}(\psi) - m_{N,t}(\psi)) + \sqrt{N}(m_{N,t}(\psi) - m_t(\psi)). \end{aligned} \tag{25}$$

We assume that, conditionally on all samples up to time $t-1$, $(x_{j,t})$ is an i.i.d. sample from $f^N_{t|t}$. Then the first term in (25) has the conditional limit distribution $\mathcal{N}(0, \sigma^2_{N,t}(\psi))$, where

$$\begin{aligned} \sigma^2_{N,t}(\psi) &= \int (\psi(x) - m_{N,t}(\psi))^2 f^N_{t|t}(x)\,d\mu(x) \\ &\approx \sigma^2_t(\psi) = \int (\psi(x) - m_t(\psi))^2 f_{t|t}(x)\,d\mu(x) \end{aligned}$$

if $f^N_{t|t}$ converges to $f_{t|t}$. By the recursions for $f_{t|t}$ and $f^N_{t|t}$, (1)–(2) and (6), respectively,

$$\sqrt{N}(m_{N,t}(\psi) - m_t(\psi)) = \sqrt{N}\left(\frac{\sum_j L_t\psi(x_{j,t-1})}{\sum_j L_t 1(x_{j,t-1})} - \frac{m_{t-1}(L_t\psi)}{m_{t-1}(L_t 1)}\right), \tag{26}$$

where

$$L_t\psi(x_{t-1}) = \int a_t(x_{t-1}, x_t) b_t(x_t, y_t)\psi(x_t)\,d\mu(x_t).$$

Asymptotic normality of the second term of (25) follows therefore from the induction assumption and the delta method.

We now state and prove a rigorous result.

THEOREM 4. *If $x \to a_t(x,\cdot)$ is continuous and if for all $t$, all $x$ and all $y$,*

$$0 < b_t(x,y) \le C(t,y) < \infty,$$

*then for all $t$, all $y_{1:t}$ and all functions $\psi$ with*

$$\sigma_t^2(\psi) = \int (\psi(x) - m_t(\psi))^2 f_{t|t}(x)\, d\mu(x) < \infty,$$

*the recursively defined asymptotic variance*

$$V_t(\psi) = \sigma_t^2(\psi) + \frac{1}{p^2(y_t|y_{1:t-1})} V_{t-1}(L_t(\psi - m_t(\psi)))$$

*is finite. Moreover, if $\sigma_s^2(\psi_s) < \infty$ for $s = 0,1,\dots,t$, then the vector $\sqrt{N} \times (M_{N,s}(\psi_s) - m_s(\psi_s))_{s=0,\dots,t}$ converges in distribution to a $\mathcal{N}(0,(V_{r,s}(\psi_r,\psi_s)))$ random vector, where*

$$V_{r,t}(\psi_r,\psi_t) = V_{r,t-1}(\psi_r, L_t(\psi_t - m_t(\psi_t)))$$

*for $r < t$ and $V_{t,t}(\psi_t,\phi_t) = (V_t(\psi_t + \phi_t) - V_t(\psi_t) - V_t(\phi_t))/2$.*

PROOF. Using the Cramér–Wold device, it is sufficient to show that

$$Z_N = \sqrt{N} \sum_{s=0}^{t} (M_{N,s}(\psi_s) - m_s(\psi_s))$$

is asymptotically centered normal with variance

$$\tau^2 = \sum_{r,s=0}^{t} V_{r,s}(\psi_r,\psi_s).$$

For $t = 0$, the theorem is trivially satisfied, and for the induction argument, we decompose $Z_N = Z_N^{(1)} + Z_N^{(2)}$, where

$$Z_N^{(1)} = \sqrt{N}(M_{N,t}(\psi_t) - m_{N,t}(\psi_t))$$

and

$$Z_N^{(2)} = \sqrt{N}(m_{N,t}(\psi_t) - m_t(\psi_t)) + \sqrt{N} \sum_{s=0}^{t-1} (M_{N,s}(\psi_s) - m_s(\psi_s)).$$

We first assume that $\psi_t$ is bounded. Denoting by $\mathcal{F}_t$ the $\sigma$-field generated by the $(x_{j,s}; 1 \le j \le N, 0 \le s \le t)$, we can write

$$\mathbf{E}[\exp(i\lambda Z_N)] = \mathbf{E}[\mathbf{E}[\exp(i\lambda Z_N^{(1)})|\mathcal{F}_{t-1}] \exp(i\lambda Z_N^{(2)})].$$

Since conditionally on $\mathcal{F}_{t-1}$ the $x_{j,t}$'s are i.i.d., we have

$$\mathbf{E}[\exp(i\lambda Z_N^{(1)})|\mathcal{F}_{t-1}] = \left(\mathbf{E}\left[\exp\left(i\frac{\lambda}{\sqrt{N}}(\psi_t(x_{1,t}) - m_{N,t}(\psi_t))\right)\Big|\mathcal{F}_{t-1}\right]\right)^N.$$

Furthermore, by a Taylor expansion of $\exp(iu)$,

$$\left|\mathbf{E}\left[\exp\left(i\frac{\lambda}{\sqrt{N}}(\psi_t(x_{1,t}) - m_{N,t}(\psi_t))\right)\Big|\mathcal{F}_{t-1}\right] - 1 + \frac{\lambda^2\sigma_{N,t}^2(\psi_t)}{2N}\right| \\ \le \frac{|\lambda|^3 \sup|\psi_t(x)|^3}{6N^{3/2}}.$$

Similarly, because $1 - u \le \exp(-u) \le 1 - u + u^2$ for all $u \ge 0$,

$$\left|1 - \frac{\lambda^2\sigma_{N,t}^2(\psi_t)}{2N} - \exp(-\lambda^2\sigma_{N,t}^2(\psi_t)/(2N))\right| \le \frac{\lambda^4 \sup|\psi_t(x)|^4}{4N^2}.$$

Because $|u^N - v^N| \le N|u - v|$ for $|u| \le 1, |v| \le 1$, we therefore obtain that, for any $\lambda$,

$$\mathbf{E}[\exp(i\lambda Z_N^{(1)})|\mathcal{F}_{t-1}] - \exp(-\lambda^2\sigma_{N,t}^2(\psi_t)/2)$$

converges to zero as $N \to \infty$ uniformly. By Theorem 1, $\|f_{t|t}^N - f_{t|t}\|_1$ converges to zero for $N \to \infty$. Because $\psi_t$ is bounded, this implies that $\sigma_{N,t}^2(\psi_t)$ converges to $\sigma_t^2(\psi_t)$. Therefore,

$$\mathbf{E}[|\exp(-\lambda^2\sigma_{N,t}^2(\psi)/2) - \exp(-\lambda^2\sigma_t^2(\psi)/2)|] \stackrel{N\to\infty}{\longrightarrow} 0.$$

We now turn to the second term, $Z_N^{(2)}$. The conditions of the theorem guarantee that

$$\begin{aligned} m_{t-1}(L_t 1) &= \iint f_{t-1|t-1}(x_{t-1}) a_t(x_{t-1}, x_t) b_t(x_t, y_t)\, d\mu(x_{t-1})\, d\mu(x_t) \\ &= p(y_t|y_{1:t-1}) \end{aligned}$$

is strictly positive, and $L_t\psi_t$ and $L_t 1$ are easily seen to be bounded if $\psi_t$ is bounded. Hence, the conditions for the delta method are satisfied, and so $\sqrt{N}(m_{N,t}(\psi_t) - m_t(\psi_t))$ is asymptotically equivalent to

$$\frac{1}{\sqrt{N} p(y_t|y_{1:t-1})}\left(\sum_j (L_t\psi_t(x_{j,t-1}) - m_{t-1}(L_t\psi_t)) \\ - m_t(\psi_t)\sum_j (L_t 1(x_{j,t-1}) - m_{t-1}(L_t 1))\right).$$

This is equal to $\sqrt{N}(M_{N,t-1}(\phi_{t-1}) - m_{t-1}(\phi_{t-1}))$, where $\phi_{t-1} = L_t(\psi_t - m_t(\psi_t))/p(y_t|y_{1:t-1})$. Hence, by the induction assumption, $\mathbf{E}[\exp(i\lambda Z_N^{(2)})]$ converges to

$$\exp\Biggl(-\frac{\lambda^2}{2}\Biggl(\sum_{r,s=0}^{t-2} V_{r,s}(\psi_r,\psi_s) + 2\sum_{s=0}^{t-2} V_{s,t-1}(\psi_s,\psi_{t-1}+\phi_{t-1}) + V_{t-1}(\psi_{t-1}+\phi_{t-1})\Biggr)\Biggr),$$

which is equal to $\exp(-\lambda^2(\tau^2 - \sigma_t^2(\psi_t))/2)$ because $V_{r,t}(\cdot,\cdot)$ is bilinear.

Taking all this together we obtain that, for bounded $\psi_t$,

$$\begin{aligned}
&|\mathbf{E}[\exp(i\lambda Z_N)] - \exp(-\lambda^2\tau^2/2)| \\
&\quad \le \mathbf{E}[|\mathbf{E}[\exp(i\lambda Z_N^{(1)})|\mathcal{F}_{t-1}] - \exp(-\lambda^2\sigma_t^2(\psi_t)/2)|] \\
&\qquad + |\mathbf{E}[\exp(i\lambda Z_N^{(2)})] - \exp(-\lambda^2(\tau^2 - \sigma_t^2(\psi_t))/2)|
\end{aligned}$$

converges to zero.

The last part of the proof deals with the case when $\psi_t$ is unbounded. We show first that $\sigma_t(\psi_t) < \infty$ implies $V_t(\psi_t) < \infty$. Again we use induction. For $t = 0$, this is clear because $\sigma_0^2(\psi) = V_0(\psi)$. For the induction step, it is sufficient to show that $\sigma_{t-1}(L_t(\psi_t - m_t(\psi_t))) < \infty$ because, by our assumptions, $p(y_t|y_{1:t-1}) > 0$. By Cauchy–Schwarz, $L_t^2\psi \le L_t(\psi^2)L_t1$, and by our assumption, $L_t1 \le C(t, y_t)$ is finite. Hence, by the definition of $L_t$ and the recursions (1)–(2),

$$\begin{aligned}
\sigma_{t-1}^2(L_t(\psi_t - m_t(\psi))) &\le m_{t-1}(L_t^2(\psi_t - m_t(\psi))) \\
&\le C(t, y_t)m_{t-1}(L_t((\psi_t - m_t(\psi))^2)) \\
&= C(t, y_t)p(y_t|y_{1:t-1})\sigma_t^2(\psi_t) < \infty.
\end{aligned}$$

For the asymptotic normality, we use a truncation argument. We set

$$\psi_{t,c}(x) = \psi_t(x)\mathbb{1}_{\{|\psi_t(x)|\le c\}}, \qquad \overline{\psi}_{t,c}(x) = \psi_t(x) - \psi_{t,c}(x).$$

Because $V_t(\psi_t) < \infty$, it follows by dominated convergence that

$$V_{r,t}(\psi_r, \psi_{t,c}) \stackrel{c\to\infty}{\longrightarrow} V_{r,t}(\psi_r, \psi_t). \tag{27}$$

Next, we are going to show that

$$\lim_{c\to\infty} \limsup_N \mathbf{P}[\sqrt{N}|M_{N,t}(\overline{\psi}_{t,c}) - m_t(\overline{\psi}_{t,c})| \ge \epsilon] = 0. \tag{28}$$

We first condition on $\mathcal{F}_{t-1}$. By Chebyshev's inequality,

$$\begin{aligned}
&\Pr[\sqrt{N}|M_{N,t}(\overline{\psi}_{t,c}) - m_t(\overline{\psi}_{t,c})| \ge \epsilon|\mathcal{F}_{t-1}] \\
&\qquad \le \mathbb{1}_{\{\sqrt{N}|m_{N,t}(\overline{\psi}_{t,c}) - m_t(\overline{\psi}_{t,c})|\ge\epsilon/2\}} + \min\biggl(1, \frac{4}{\epsilon^2}m_{N,t}(\overline{\psi}_{t,c}^2)\biggr).
\end{aligned}$$

We therefore have to study the expectations of the two terms on the right. By (26),

$$\sqrt{N}(m_{N,t}(\overline{\psi}_{t,c}) - m_t(\overline{\psi}_{t,c})) = \sqrt{N}\left(\frac{\sum_j L_t\overline{\psi}_{t,c}(x_{j,t-1})}{\sum_j L_t 1(x_{j,t-1})} - \frac{m_{t-1}(L_t\overline{\psi}_{t,c})}{m_{t-1}(L_t 1)}\right),$$

which by the induction assumption is asymptotically $\mathcal{N}(0, V_{t-1}(L_t(\overline{\psi}_{t,c} - m_t(\overline{\psi}_{t,c}))))$-distributed. For $c \to \infty$, this variance goes to zero, implying the desired behavior of the first term. By the recursion for $f^N_{t|t}$,

$$m_{N,t}(\overline{\psi}^2_{t,c}) = \frac{\sum_j L_t\overline{\psi}^2_{t,c}(x_{j,t-1})}{\sum_j L_t 1(x_{j,t-1})},$$

which, by the induction assumption, converges in probability to $\int \overline{\psi}^2_{t,c}(x) f_{t|t}(x)\, d\mu(x)$. Hence, by dominated convergence the second term also has the desired behavior, and, thus, (28) follows.

Now we have all the ingredients to complete the proof. We write

$$Z_{N,c} = \sqrt{N}\sum_{s=0}^{t-1}(M_{N,s}(\psi_s) - m_s(\psi_s)) + (M_{N,t}(\psi_{t,c}) - m_t(\psi_{t,c}))$$

and $\tau_c^2$ for the asymptotic variance of $Z_{N,c}$. Then

$$\begin{aligned}
&|\mathbf{E}[\exp(i\lambda Z_N)] - \exp(-\lambda^2\tau^2/2)| \\
&\quad \le |\mathbf{E}[\exp(i\lambda Z_{N,c})] - \exp(-\lambda^2\tau_c^2/2)| \\
&\qquad + |\exp(-\lambda^2\tau_c^2/2) - \exp(-\lambda^2\tau^2/2)| \\
&\qquad + \mathbf{E}[|\exp(i\lambda\sqrt{N}(M_{N,t}(\psi_{t,c}) - m_t(\psi_{t,c}) - M_{N,t}(\psi) + m_t(\psi))) - 1|].
\end{aligned}$$

By (27), the second term is arbitrarily small if $c$ is large enough. Using $|\exp(iu) - 1| \le \min(2, |u|)$ and (28), the same thing holds also for the last term, uniformly in $N$. Finally, the first term goes to zero for any fixed $c$ as $N \to \infty$. □

4.3.1. *The asymptotic variance.* Similarly as in the case of convergence of $f^N_{t|t}$, one would like to know whether the asymptotic variances $V_t(\psi)$ stay bounded as $t$ increases. Using ideas from [3], we show that this is the case if $\psi$ is bounded and the condition (19) is satisfied. Because $m_{t-1}(L_t\psi) = m_t(\psi)p(y_t|y_{1:t-1})$, we have

$$m_{t-1}(L_t(\psi - m_t(\psi))) = 0.$$

Hence, by iterating the recursive definition of $V_t(\psi)$, we obtain

$$V_t(\psi) = \sigma_t^2(\psi) + \sum_{s=1}^{t}\frac{\sigma^2_{s-1}(L_{s:t}(\psi - m_t(\psi)))}{p^2(y_{s:t}|y_{1:s-1})}, \tag{29}$$

where

$$L_{s:t}\psi(x_{s-1}) = \int \psi(x_t) \prod_{r=s}^{t} a_r(x_{r-1}, x_r) b_r(x_r, y_r)\, d\mu(x_r)$$
$$= \mathbf{E}[\psi(X_t)|x_{s-1}, y_{s:t}] p(y_{s:t}|x_{s-1}).$$

Here, the expectation is with respect to the state space model and not with respect to the random sampling in the Monte Carlo filter. Thus,

$$L_{s:t}(\psi - m_t(\psi)) = (\mathbf{E}[\psi(X_t)|x_{s-1}, y_{s:t}] - \mathbf{E}[\psi(X_t)|y_{1:t}]) p(y_{s:t}|x_{s-1}).$$

Because $p(y_{s:t}|y_{1:s-1}) = \int p(y_{s:t}|x_{s-1}) f_{s-1|s-1}(x_{s-1}|y_{1:s-1})\, d\mu(x_{s-1})$, it follows from (22) that

$$\frac{p(y_{s:t}|x_{s-1})}{p(y_{s:t}|y_{1:s-1})} \leq \frac{1}{\gamma_a}.$$

Moreover, condition (19) implies uniform contractivity of the conditional chain given $y_s^{t-1}$; compare Lemma 10. Hence, under condition (19) we have

$$|\mathbf{E}[\psi(X_t)|x_{s-1}, y_{s:t}] - \mathbf{E}[\psi(X_t)|y_{1:t}]| \leq \left(\sup_x \psi(x) - \inf_x \psi(x)\right)(1-\gamma_a)^{t-s+1}$$

and, therefore,

$$V_t(\psi) \leq \gamma_a^{-3} \left(\sup_x \psi(x) - \inf_x \psi(x)\right)^2. \tag{30}$$

So far, we have dealt with the case where $(x_{j,t})$ is an i.i.d. sample from $f_{t|t}^N$, usually generated by an accept–reject method. For sampling importance resampling, asymptotic normality can be proved by a similar recursive argument; see [3]. However, the formula for the variance $V_t$ changes slightly. Random resampling leads to the recursion

$$V_t(\psi) = \frac{1}{p^2(y_t|y_{1:t-1})} V_{t-1}(L_t(\psi - m_t(\psi))) + \sigma_t^2(\psi)$$
$$+ \frac{1}{p^2(y_t|y_{1:t-1})} m_{t-1}(L_t(b_t(\cdot y_t)(\psi - m_t(\psi))^2) - L_t^2(\psi - m_t(\psi))).$$

(The second term comes from the resampling step and the third from the reweighting.) Using again $m_{t-1}(L_t\psi) = m_t(\psi) p(y_t|y_{1:t-1})$, we obtain

$$V_t(\psi) = \sigma_t^2(\psi) + \frac{V_{t-1}(L_t(\psi - m_t(\psi))) - \sigma_{t-1}^2(L_t(\psi - m_t(\psi)))}{p^2(y_t|y_{1:t-1})}$$
$$+ \frac{m_t(b_t(\psi - m_t(\psi))^2)}{p(y_t|y_{1:t-1})} \tag{31}$$
$$= \sigma_t^2(\psi) + \sum_{s=1}^{t} \frac{m_s(b_s L_{s+1:t}^2(\psi - m_t(\psi)))}{p(y_s|y_{1:s-1}) p^2(y_{s+1:t}|y_{1:s})}$$

[we set $L_{t+1:t}(\psi) = \psi$]. Using Cauchy–Schwarz, one can show that each summand in (31) is always greater than or equal to the corresponding term in (29) and, thus, the additional effort of generating an i.i.d. sampling reduces the variance.

Because of the slightly different form of the asymptotic variance, one also needs additional conditions in order that $V_t(\psi)$ in (31) remain bounded uniformly in $t$. Using the previous bound for $(L_{(s+1):t}(\psi - m_t(\psi)))/p(y_{s+1:t}|y_{1:s})$, one needs, in addition, a bound for

$$\frac{m_s(b_s)}{p(y_s|y_{1:s-1})} = \frac{\int f_{s-1|s-1}(x_{s-1}) a_s(x_{s-1}, x_s) b_s^2(x_s, y_s)\, d\mu(x_{s-1})\, d\mu(x_s)}{(\int f_{s-1|s-1}(x_{s-1}) a_s(x_{s-1}, x_s) b_s(x_s, y_s)\, d\mu(x_{s-1})\, d\mu(x_s))^2}.$$

Obviously, this is bounded uniformly in $s$ and $y$ under the condition (21). Using assumption (19), we can replace (21) by a slightly stronger version of (20), namely, that

$$\int h(x_s) b_s^2(x_s, y_s)\, d\mu(x_s) < \infty \tag{32}$$

for all $s$ and all $y_s$. However, the bound for $V_t$ then depends on $y_{1:t}$.

4.4. *High rate sampling.* So far, we have worked with a fixed sampling rate which we set equal to one for simplicity. Alternatively, we can consider what happens when the sampling rate converges to zero. We discuss this case briefly in this last section. So we let $(X_t)$ be a Markov process in continuous time, and we assume, for simplicity, that it is time homogeneous with transition kernels $\mathbf{P}[X_{t+s} \in dx | X_t = x'] = a(s, x', x)\, d\mu(x)$. We consider the sampling rate $\delta = 1/m$ with $m \in \mathbf{N}$, and we assume that, for a given $m$, we have conditionally independent observations $Y_{j\delta}, j = 1, 2, \ldots,$ such that $Y_{j\delta}$ depends only on $X_{j\delta}$.

In the previous two subsections we showed how the strong condition (19) allows one to obtain convergence results that are uniform in $t$ and require essentially no conditions on the observation densities. Unfortunately, this strategy breaks down in the high rate sampling limit. In continuous time, the analogue of (19) is

$$c_a(t) h(x) \le a(t, x', x) \le C_a(t) h(x) \tag{33}$$

for some fixed $h$ and all $t, x, x'$. It is easy to see that if the lower bound $c_a(t)$ is of larger order than $t$ as $t \to 0$, then $\|a_t(x, \cdot) - a_t(x', \cdot)\|_1 = 0$ for all $t > 0$. Hence, except for degenerate cases, the crucial quantity $\gamma_a$ diverges at least like $\delta^{-1}$ for $\delta \to 0$. Moreover, the continuity module of $x \to a(\delta, x, \cdot)$ which is used in Lemmas 8 and 9 also diverges. Because the asymptotic variance $V_t(\psi)$ in (29) is exact and does not depend on Lemmas 8 and 9, it is slightly easier to study the behavior of $V_t(\psi)$ as $\delta \to 0$, and we concentrate on this.

Even this is not trivial. The simplest case occurs if the state space is finite and all jump probabilities are positive. Then it is easily seen that (33) holds with $c_a(t) = c_a t$, $c_a > 0$ some constant, and $C_a(t) \equiv C_a$. Inserting this into the bound (30) for the asymptotic variance, we obtain an upper bound for the asymptotic variance of the order $m^3$ which is not satisfactory. Of course, our bounds are presumably not sharp, but it is not obvious how to improve them in general. We believe that the behavior in the high rate sampling case depends on the properties of the observation process. If we have a fixed observation density and we increase the sampling rate, we accumulate more and more information about the state process in any fixed interval, and the filter distribution will converge to a point mass except near the times where a jump occurs. With high rate sampling, it is somewhat more natural to let the information that is carried by a single observation decrease with the sampling rate. Then we need additional superscripts for the observations and their densities. The standard example is

$$Y_{j\delta}^{(m)} | X_{j\delta} = x \sim \mathcal{N}(\delta g(x), \delta\sigma^2), \tag{34}$$

and we will study this case. Then the partial sum process

$$\eta_t^{(m)} = \sum_{j\delta \le t} Y_{j\delta}^{(m)} \tag{35}$$

converges for $m \to \infty$ in distribution to the process $\eta_t = \int_0^t g(x_s)\, ds + \sigma W_t$.

If $t$ is fixed and the sampling rate increases, the formula for $V_t$ contains $O(m)$ summands. Moreover, with (34), the filter distributions are not degenerate and the function $x_{s-\delta} \to \mathbf{E}[\psi(X_t) | x_{s-\delta}, y_{s:t}]$ does not converge to a constant. Hence, we expect that, for fixed $t, V_t(\psi)$ is of the order $m$. This is not surprising: At each time step $\delta$ we take a new sample even though the filter distribution changes very little. The sampling errors accumulate because the filter does not forget its initial condition over a finite time interval. In this setup, it is much better to use sequential importance sampling, that is, to carry the weights forward by multiplication instead of resampling at each time step. We thus generate a sample $(x_{j,k\delta}, k = 0, 1, \dots)$ from our model of the state process and compute the weights $\Lambda_{1:k}^{(m)}(x_{j,\delta:k\delta})$ sequentially, where

$$\Lambda_{i:k}^{(m)}(x_{i\delta:k\delta}) = \prod_{\ell=i}^{k} b^{(m)}(x_{\ell\delta}, y_{\ell\delta}^{(m)})$$

is the likelihood. Then

$$M_{N,k\delta}(\psi) = \frac{\sum_{j=1}^N \Lambda_{1:k}^{(m)}(x_{j,\delta:k\delta})\psi(x_{j,k\delta})}{\sum_{j=1}^N \Lambda_{1:k}^{(m)}(x_{j,\delta:k\delta})}$$

is an asymptotically normal estimator of $m_{k\delta}(\psi)$ with asymptotic variance

$$V_{k\delta}(\psi)=\frac{\mathbf{E}_X[(\psi(X_{k\delta})-m_{k\delta}(\psi))^2\Lambda_{1:k}^{(m)2}(X_{\delta:k\delta})]}{(\mathbf{E}_X[\Lambda_{1:k}^{(m)}(X_{\delta:k\delta})])^2}$$

($\mathbf{E}_X$ indicates that the expectation is only with respect to the state variables, the observations are considered to be fixed). We show first that, for a fixed time $k\delta$, this variance remains bounded as the sampling rate increases to infinity. We assume the state process to be a diffusion,

$$dX_t=f(X_t)\,dt+\tilde\sigma(X_t)\,dB_t, \tag{36}$$

where $(B_t)$ is a Brownian motion.

THEOREM 5. *Consider the state space model $(X_{j\delta},Y_{j\delta}^{(m)})$ defined by (36) and (34), where $f$, $\tilde\sigma$ and $g$, together with their first and second derivatives, are all bounded. Assume, moreover, that the partial sum process (35) converges in the sup-norm to a continuous function $\eta$. Then for $\delta\to 0$ and $k\delta\to t$, $t$ fixed, the asymptotic variance $V_{k\delta}(\psi)$ of the importance sampling estimator is bounded by $\sup_x\psi^2(x)$ times a constant that depends only on $t$ and the supremum of $|\eta|$ on $[0,t]$.*

PROOF. In order to simplify the notation we assume that $\sigma=1$. Moreover, we put

$$\begin{aligned}\overline{b}(x,y)&=\exp(-\delta g(x)^2/2+g(x)y)\\&=b(x,y)\sqrt{(2\pi)/m}\exp(my^2/2).\end{aligned}\tag{37}$$

Then we can replace the observation density $b$ in the likelihood $\Lambda^{(m)}$ by $\overline{b}$.

Summation by parts gives

$$\begin{aligned}\Lambda_{1:k}^{(m)}(x_{\delta:k\delta})=\exp\Bigg(&-\delta/2\sum_{\ell=1}^{k}g(x_{\ell\delta})^2\\&-\sum_{\ell=1}^{k-1}(g(x_{(\ell+1)\delta})-g(x_{\ell\delta}))\eta_{\ell\delta}^{(m)}+g(x_{k\delta})\eta_{k\delta}^{(m)}\Bigg).\end{aligned}$$

The Itô formula implies

$$\begin{aligned}&g(x_{(\ell+1)\delta})-g(x_{\ell\delta})\\&\quad=\int_{\ell\delta}^{(\ell+1)\delta}Dg(x_t)\tilde\sigma(x_t)\,dB_t\\&\qquad+\int_{\ell\delta}^{(\ell+1)\delta}(Dg(x_t)f(x_t)+\tfrac12\tilde\sigma(x_t)^TD^2g(x_t)\tilde\sigma(x_t))\,dt.\end{aligned}$$

Because $f$, $\tilde{\sigma}$, $g$, $Dg$ and $D^2g$ are all bounded, it follows that

$$\log \Lambda_{1:k}^{(m)}(x_{\delta:k\delta}) + \int_0^{k\delta} \eta_t^{(m)} Dg(x_t)\tilde{\sigma}(x_t)\,dB_t$$

is bounded above and below by constants that depend only on $k\delta$ and $\sup\{|\eta_t^{(m)}|; 0 \le t \le k\delta\}$. Finally, again by the Itô formula,

$$\begin{aligned}
&\mathbf{E}\biggl[\exp\biggl(\int_0^{k\delta} \eta_t^{(m)} Dg(X_t)\tilde{\sigma}(X_t)\,dB_t\biggr)\biggr] \\
&\qquad = \mathbf{E}\biggl[\exp\biggl(-\tfrac{1}{2}\int_0^{k\delta} \eta_t^{(m)2}(Dg(X_t)\tilde{\sigma}(X_t))^2\,dt\biggr)\biggr],
\end{aligned}$$

which is again bounded above and below by constants that depend only on $k\delta$ and the supremum of $|\eta_t^{(m)}|$. Using this, the rest of the proof is obvious. □

Finally, one can look at the case where both the sampling rate $m$ and the time index $t$ of the filtering distribution increase. We show that in this situation the asymptotic variance remains bounded if we resample at fixed time intervals which for simplicity we take equal to one. Accept–reject methods at a fixed rate cannot be used in this case because the supremum of $\Lambda_{(s-1)m+1:sm}^{(m)}(x_{s-1+\delta:s})$ diverges as $m$ increases.

By similar arguments as before, the asymptotic variance at time $t \in \mathbb{N}$ of this version of the particle filter is

$$\begin{aligned}
V_t(\psi) &= \sigma_t^2(\psi) \\
&\quad + \sum_{s=1}^{t}[m_{s-1}(\mathbf{E}_X[\Lambda_{(s-1)m+1:sm}^{(m)2}(X_{s-1+\delta:s}) \\
&\qquad\qquad \times L_{s+\delta:t}^{(m)2}(\psi - m_t(\psi))(X_s)|x_{s-1}]) \\
&\qquad\qquad\qquad \times [p(y_{s-1+\delta:t}^{(m)}|y_{\delta:s-1}^{(m)})]^{-2}].
\end{aligned} \tag{38}$$

As above, we replace the observation density $b$ in the likelihood $\Lambda^{(m)}$ by $\overline{b}$ since this has no effect on the right-hand side of (38). We define

$$\begin{aligned}
J_{j:k}^{(m)}(x_{j\delta}, x_{k\delta}) = \int &\Lambda_{j+1:k}^{(m)}(x_{(j+1)\delta:k\delta})a(\delta, x_{(k-1)\delta}, x_{k\delta}) \\
&\times \prod_{i=j+1}^{k-1} a(\delta, x_{(i-1)\delta}, x_{i\delta})\,d\mu(x_{i\delta}).
\end{aligned}$$

We then assume that there exist a probability density $h$ and two functions $c$ and $C$ $\mathbb{R}_+^2 \to \mathbb{R}_+$ such that, for all $x_{j\delta}$,

$$c((j-k)\delta, M^{(m)}(j\delta, k\delta)) \leq \frac{J_{j:k}^{(m)}(x_{j\delta}, x_{k\delta})}{h(x_{k\delta})} \leq C((j-k)\delta, M^{(m)}(j\delta, k\delta)), \tag{39}$$

where

$$M^{(m)}(s,t) = \sup_{s<u\leq t} |\eta_u^{(m)} - \eta_s^{(m)}|.$$

It follows from the proof of Theorem 3 in [1] that the assumption (39) is satisfied in the case where the state process is a diffusion on a compact Riemannian manifold with strictly elliptic generator.

Because for $j < \ell < k$

$$J_{j:k}^{(m)}(x_{j\delta}, x_{k\delta}) = \int J_{j:\ell}^{(m)}(x_{j\delta}, x_{\ell\delta}) J_{\ell:k}^{(m)}(x_{\ell\delta}, x_{k\delta})\, d\mu(x_{\ell\delta}),$$

assumption (39) implies that, for $k - j \geq m$,

$$\frac{p(y_{(j+1)\delta:k\delta}^{(m)} | x_{j\delta})}{p(y_{(j+1)\delta:k\delta}^{(m)} | x'_{j\delta})} = \frac{\int J_{j:k}^{(m)}(x_{j\delta}, x_{k\delta})\, d\mu(x_{k\delta})}{\int J_{j:k}^{(m)}(x'_{j\delta}, x_{k\delta})\, d\mu(x_{k\delta})} \leq \frac{C(1, M^{(m)}(j\delta, j\delta+1))}{c(1, M^{(m)}(j\delta, j\delta+1))}.$$

Similarly,

$$p(x_{j\delta+1} | x_{j\delta}, y_{(j+1)\delta:k\delta}^{(m)}) = \frac{\int J_{j:j+m}^{(m)}(x_{j\delta}, x_{j\delta+1}) J_{j+m:k}^{(m)}(x_{j\delta+1}, x_{k\delta})\, d\mu(x_{k\delta})}{\int J_{j:k}^{(m)}(x_{j\delta}, x_{k\delta})\, d\mu(x_{k\delta})}$$

is bounded below by

$$\frac{c(1, M^{(m)}(j\delta, j\delta+1))}{C(1, M^{(m)}(j\delta, j\delta+1))} \frac{h(x_{j\delta+1}) \int J_{j+m:k}^{(m)}(x_{j\delta+1}, x_{k\delta})\, d\mu(x_{k\delta})}{\int h(x_{j\delta+1}) J_{j+m:k}^{(m)}(x_{j\delta+1}, x_{k\delta})\, d\mu(x_{j\delta+1})\, d\mu(x_{k\delta})}.$$

The second ratio on the right-hand side is a probability density which does not depend on $x_{j\delta}$, and, thus, one minus the first ratio is a bound for the contraction rate of the conditional chain. Therefore, we have

$$\frac{|L_{s+\delta:t}^{(m)}(\psi - m_t(\psi))(x_s)|}{p(y_{s+\delta:t}^{(m)} | y_{\delta:s})} \leq \frac{\sup_x \psi(x) - \inf_x \psi(x)}{\gamma(s)} \prod_{r=s}^{t-1} (1 - \gamma(r)),$$

where

$$\gamma(r) = \frac{c(1, M^{(m)}(r, r+1))}{C(1, M^{(m)}(r, r+1))}.$$

In contrast to the case with fixed sampling rate, these coefficients depend on the observations, that is, they are random. For nonrandom bounds that hold with high probability, we would have to assume that the limiting process $\eta$ has stationary increments. Finally, we can bound

$$\frac{m_{s-1}(\mathbf{E}_X[\Lambda^{(m)2}_{(s-1)m+1\,:\,sm}(X_{s-1+\delta\,:\,s})|x_{s-1}])}{p^2(y^{(m)}_{s-1+\delta\,:\,s}|y^{(m)}_{\delta\,:\,s-1})}$$

by similar arguments as before.

**Acknowledgments.** I am grateful to Neil Shephard, Eric Moulines and two referees for helpful comments on earlier versions of this paper. In particular, I thank Eric Moulines for pointing out to me that Theorem 3 holds.

Seminar für Statistik
ETH Zentrum
CH-8092 Zürich
Switzerland
e-mail: kuensch@stat.math.ethz.ch